\documentclass{amsart}
\usepackage{amssymb,amscd,pstricks,pst-node,pst-tree}
\usepackage[stdtext]{youngtab}

\numberwithin{equation}{section}

\newtheorem{prop}{Proposition}
\newtheorem{theorem}[prop]{Theorem}

\newtheorem{conjecture}[prop]{Conjecture}
\newtheorem{defconj}[prop]{Definition-Conjecture}

\theoremstyle{definition}
\newtheorem{definition}[prop]{Definition}
\newtheorem{example}[prop]{Example}

\numberwithin{prop}{section}

\newcommand{\0}{\circ}
\newcommand{\alt}{\tilde{\alpha}}
\newcommand{\aut}{\sigma}
\newcommand{\bij}{\iota}
\newcommand{\cc}{\mathrm{cc}}
\newcommand{\dd}{\mathrm{d}}
\newcommand{\dom}{\trianglerighteq}
\newcommand{\emb}{\Psi}
\newcommand{\eh}{\widehat{e}}
\newcommand{\es}{\varnothing}
\newcommand{\fh}{\widehat{f}}
\newcommand{\geh}{\mathfrak{g}}
\newcommand{\gehb}{\overline{\geh}}

\newcommand{\inner}[2]{\langle #1\,,\,#2\rangle}
\newcommand{\la}{\lambda}
\newcommand{\lev}{\mathrm{lev}}
\newcommand{\mh}{\widehat{m}}
\newcommand{\mult}{\gamma}
\newcommand{\nat}{u^\natural}
\newcommand{\nh}{\widehat{\nu}}
\newcommand{\ph}{\widehat{p}}
\newcommand{\qbin}[2]{\genfrac{[}{]}{0pt}{}{#1}{#2}}
\newcommand{\uu}{u}
\newcommand{\ve}[1]{\varepsilon_{#1}}
\newcommand{\viso}{\emb}
\newcommand{\vp}[1]{\varphi_{#1}}
\newcommand{\wt}{\mathrm{wt}}
\newcommand{\xb}{\bar{x}}
\newcommand{\yv}{y^\vee}

\newcommand{\Cv}{C^v}
\newcommand{\Dv}{D^v}
\newcommand{\Hv}{H^v}
\newcommand{\I}{I}
\newcommand{\J}{\overline{\I}}
\newcommand{\Jh}{\widehat{J}}
\newcommand{\La}{\Lambda}
\newcommand{\Lab}{\overline{\La}}
\newcommand{\Lh}{\widehat{L}}
\newcommand{\Path}{\mathcal{P}}
\newcommand{\Pfin}{\overline{P}}
\newcommand{\Qfin}{\overline{Q}}
\newcommand{\R}{R}
\newcommand{\Rh}{\widehat{R}}
\newcommand{\RC}{\mathrm{RC}}
\newcommand{\RCv}{\RC^v}
\newcommand{\Rv}{R^v}
\newcommand{\Th}{\widehat{T}}
\newcommand{\V}{V}
\newcommand{\Vh}{\widehat{\V}}
\newcommand{\Wb}{\overline{W}}
\newcommand{\X}{\mathcal{X}}
\newcommand{\Xv}{X^v}
\newcommand{\Y}{\mathcal{Y}}
\newcommand{\Yv}{\Y^\vee}
\newcommand{\Z}{\mathbb{Z}}

\begin{document}

\title{Virtual crystals and Kleber's algorithm}

\author[M.~Okado]{Masato Okado}
\address{Department of Informatics and Mathematical Science,
Graduate School of Engineering Science, Osaka University,
Toyonaka, Osaka 560-8531, Japan}
\email{okado@sigmath.es.osaka-u.ac.jp}

\author[A.~Schilling]{Anne Schilling}
\address{Department of Mathematics, University of California, One Shields
Avenue, Davis, CA 95616-8633, U.S.A.}
\email{anne@math.ucdavis.edu}

\author[M.~Shimozono]{Mark Shimozono}
\address{Department of Mathematics, 460 McBryde Hall, Virginia Tech,
Blacksburg, VA 24061-0123, U.S.A}
\email{mshimo@math.vt.edu}



\begin{abstract}
Kirillov and Reshetikhin conjectured what is
now known as the fermionic formula for the decomposition
of tensor products of certain finite dimensional
modules over quantum affine algebras.
This formula can also be extended to
the case of $q$-deformations of tensor product
multiplicities as recently conjectured by
Hatayama et al.. In its original formulation it is difficult
to compute the fermionic formula efficiently. Kleber found an
algorithm for the simply-laced algebras which overcomes
this problem. We present a method which reduces all
other cases to the simply-laced case using
embeddings of affine algebras. This is the fermionic
analogue of the virtual crystal construction by the authors,
which is the realization of crystal graphs for arbitrary quantum
affine algebras in terms of those of simply-laced type.
\end{abstract}

\maketitle

\section{Introduction}

In 1987 Kirillov and Reshetikhin \cite{KR:1987} conjectured a
formula, now known as the fermionic formula, for the decomposition
of tensor products of certain finite dimensional representations
over an untwisted quantum affine algebra $U_q(\geh)$ into its
$U_q(\gehb)$ components, where $\gehb$ is the simple Lie algebra
associated with the affine Kac-Moody algebra $\geh$. The
conjecture is motivated by Bethe Ansatz studies. Recently,
conjectures for fermionic formulas have been extended to
$q$-deformations of tensor product multiplicities
\cite{HKOTY:1999,HKOTT:2002}. In type $A_n^{(1)}$ this $q$-tensor
multiplicity formula appeared in \cite{KR:1988}. For a single tensor factor,
the fermionic formula gives the $\gehb$-isotypical components
of a $U_q(\geh)$-module associated with a multiple of a fundamental
weight. This conjecture was proven by Chari \cite{C:2001} in a number of cases.
Recently, Nakajima \cite{N:2002} showed in the simply-laced case that
the characters of such modules satisfy a certain system of algebraic relations
($Q$-system). Combining the result of \cite{HKOTY:1999}, his result completes
the proof of a ``weak" version of the $q=1$ fermionic formula in this case.

The term fermionic formula was coined by the Stony Brook
group \cite{KKMM:1993,KKMM:1993a}, who interpreted fermionic-type
formulas for characters and branching functions of conformal
field theory models as partition functions of quasiparticle
systems with ``fractional'' statistics obeying Pauli's exclusion
principle.

Fermionic formulas are $q$-polynomials or $q$-series expressed
as certain sums of products of $q$-binomial coefficients
\begin{equation*}
\sum_{\{m_i^{(a)}\}} q^{\cc(\{m_i^{(a)}\})} \prod_{i,a}
 \qbin{m_i^{(a)}+p_i^{(a)}}{m_i^{(a)}},
\end{equation*}
where $\qbin{m+p}{m}=(q)_{m+p}/(q)_m(q)_p$ is the $q$-binomial
coefficient with $(q)_m=\prod_{i=1}^m (1-q^i)$,
$\cc(\{m_i^{(a)}\})$ is some function of the summation variables
$m_i^{(a)}$ and $p_i^{(a)}$ is the vacancy number (see
\eqref{eq:vac}). The summation variables are subject to
constraints \eqref{eq:config}. Those sets $\{m_i^{(a)}\}$
satisfying \eqref{eq:config} are called admissible configurations. {}From
\eqref{eq:config} alone, it is computationally difficult to find
the admissible configurations, making the evaluation of the
fermionic formula intractable. For simply-laced algebras $\geh$,
Kleber \cite{Kl:1997,Kl:1998} has given an efficient algorithm to
determine the admissible configurations $\{m_i^{(a)}\}$. This
algorithm generates a rooted tree with nodes labelled by dominant
integral weights such that the tree nodes are in bijection with
the admissible configurations. For non-simply laced algebras, the
algorithm fails: some admissible nodes cannot be reached.

One of our goals in this paper
is to modify Kleber's algorithm to work in all types.
This is accomplished by using the well-known natural embeddings
of any affine algebra into another of simply-laced type \cite{JM:1985}:
\begin{equation}\label{eq:embed}
\begin{array}{cll}
C_n^{(1)}, A_{2n}^{(2)}, A_{2n}^{(2)\dagger},D_{n+1}^{(2)}
&\hookrightarrow
&A_{2n-1}^{(1)}\\
A_{2n-1}^{(2)}, B_n^{(1)} &\hookrightarrow &D_{n+1}^{(1)}\\
E_6^{(2)}, F_4^{(1)} &\hookrightarrow &E_6^{(1)}\\
D_4^{(3)}, G_2^{(1)} &\hookrightarrow &D_4^{(1)}.
\end{array}
\end{equation}
It is not hard to express the fermionic formula of the smaller
algebra in terms of the larger; we call this the virtual fermionic
formula. Our algorithm is an adaptation of Kleber's algorithm
in the simply-laced affine algebra, which trims the
tree so as not to generate nodes that cannot contribute to the virtual
fermionic formula. This algorithm succeeds by using some
nodes in the larger weight lattice that do not correspond to
weights in the embedded weight lattice.

Fermionic formulas denoted $M$ have crystal counterparts. Crystal bases
were introduced by Kashiwara \cite{K:1991} and are bases of
$U_q(\geh)$-modules in the limit $q\to 0$. Let us denote
the one-dimensional configuration sums, which are generating
functions of highest weight elements in tensor products of
finite dimensional crystals with energy statistics, by $X$.
It was conjectured in \cite{HKOTY:1999,HKOTT:2002} that $X=M$.

In light of the embeddings of affine algebras \eqref{eq:embed},
one might hope that such embeddings also exist for the quantized
algebras. Unfortunately they do not. However we assert that such
embeddings exist for all finite-dimensional affine crystals, and
give a construction for them in terms of crystals of simply-laced type.
A virtual crystal is such a realization of a crystal inside another
of possibly different type. Perhaps the first instance of a virtual
crystal is Kashiwara's embedding of a crystal of highest
weight $\lambda$, into that of highest weight $k\lambda$
where $k$ is a positive integer \cite{K:1996}.
Extending Baker's work \cite{B:2000},
in \cite{OSS:2001} we conjectured that finite dimensional crystals
of type $C_n^{(1)}$, $A_{2n}^{(2)}$, and $D_{n+1}^{(2)}$ can
be realized in terms of crystals of type $A_{2n-1}^{(1)}$.
We proved this for crystals associated with single
columns (i.e. fundamental weights).

In this paper we establish the correctness of the virtual crystal
approach for crystals associated with single rows
(that is, multiples of the first
fundamental weight) for the two infinite families of embeddings.

The paper is organized as follows. In Section \ref{sec:crystals}
we review the essentials of crystal theory. Virtual crystals
are introduced in Section \ref{sec:virtual} and the characterization
and validity of virtual crystals associated with single rows is proven.
Sections \ref{sec:fermionic} and \ref{sec:kleber} review the fermionic
formulas conjectured in \cite{HKOTY:1999,HKOTT:2002} and the Kleber
algorithm, respectively, and describe their virtual counterparts.

\subsection*{Acknowledgements} Most of this work was carried out as part of
the Research in Pairs program of the Mathematisches Forschungsinstitut
Oberwohlfach in August 2002. AS and MS would like to thank the institute
for the ideal working conditions during their stay. AS also thanks the
University of Wuppertal and the Max-Planck-Institut f\"ur
Mathematik in Bonn for hospitality, where this work was completed.
MO was partially supported by Grant-in-Aid for Scientific Research
(No.14540026), JSPS. AS was partially supported by the Humboldt foundation and
NSF grant DMS-0200774. MS was partially supported by NSF grant DMS-0100918.

\section{Crystals}
\label{sec:crystals}

\subsection{Affine algebras}
We adopt the notation of \cite{HKOTT:2002}. Let
$\geh$ be a Kac-Moody Lie algebra of affine type
$X^{(r)}_N$, that is, one of the types $A^{(1)}_n (n \ge 1)$,
$B^{(1)}_n (n \ge 3)$, $C^{(1)}_n (n \ge 2)$, $D^{(1)}_n (n \ge 4)$,
$E_n^{(1)} (n=6,7,8)$, $F_4^{(1)}$, $G_2^{(1)}$,
$A^{(2)}_{2n} (n\ge1)$, $A^{(2)\dagger}_{2n} (n\ge1)$,
$A^{(2)}_{2n-1} (n\ge2)$, $D^{(2)}_{n+1}(n\ge2)$,
$E_6^{(2)}$ or $D_4^{(3)}$.
The Dynkin diagram of $\geh = X^{(r)}_N$ is depicted in Figure
\ref{fig:Dynkin} (Table Aff 1-3 in \cite{Kac}). Its nodes are
labelled by the set $\I=\{0,1,2\dotsc,n\}$.
Let $\J=\I\backslash\{0\}$.

{\unitlength=.8pt
\begin{figure}
\begin{tabular}[t]{rl}
$A_1^{(1)}$:&
\begin{picture}(26,20)(-5,-5)
\multiput( 0,0)(20,0){2}{\circle{6}}
\multiput(2.85,-1)(0,2){2}{\line(1,0){14.3}}
\put(0,-5){\makebox(0,0)[t]{$0$}}
\put(20,-5){\makebox(0,0)[t]{$1$}}
\put( 6, 0){\makebox(0,0){$<$}}
\put(14, 0){\makebox(0,0){$>$}}
\end{picture}
\\
&
\\
\begin{minipage}[b]{4em}
\begin{flushright}
$A_n^{(1)}$:\\$(n \ge 2)$
\end{flushright}
\end{minipage}&
\begin{picture}(106,40)(-5,-5)
\multiput( 0,0)(20,0){2}{\circle{6}}
\multiput(80,0)(20,0){2}{\circle{6}}
\put(50,20){\circle{6}}
\multiput( 3,0)(20,0){2}{\line(1,0){14}}
\multiput(63,0)(20,0){2}{\line(1,0){14}}
\multiput(39,0)(4,0){6}{\line(1,0){2}}
\put(2.78543,1.1142){\line(5,2){44.429}}
\put(52.78543,18.8858){\line(5,-2){44.429}}
\put(0,-5){\makebox(0,0)[t]{$1$}}
\put(20,-5){\makebox(0,0)[t]{$2$}}
\put(80,-5){\makebox(0,0)[t]{$n\!\! -\!\! 1$}}
\put(100,-5){\makebox(0,0)[t]{$n$}}
\put(55,20){\makebox(0,0)[lb]{$0$}}
\end{picture}
\\
&
\\
\begin{minipage}[b]{4em}
\begin{flushright}
$B_n^{(1)}$:\\$(n \ge 3)$
\end{flushright}
\end{minipage}&
\begin{picture}(126,40)(-5,-5)
\multiput(0,0)(20,0){3}{\circle{6}}
\multiput(100,0)(20,0){2}{\circle{6}}
\put(20,20){\circle{6}}
\multiput( 3,0)(20,0){3}{\line(1,0){14}}
\multiput(83,0)(20,0){1}{\line(1,0){14}}
\put(20,3){\line(0,1){14}}
\multiput(102.85,-1)(0,2){2}{\line(1,0){14.3}} 
\multiput(59,0)(4,0){6}{\line(1,0){2}} 
\put(110,0){\makebox(0,0){$>$}}
\put(0,-5){\makebox(0,0)[t]{$1$}}
\put(20,-5){\makebox(0,0)[t]{$2$}}
\put(40,-5){\makebox(0,0)[t]{$3$}}
\put(100,-5){\makebox(0,0)[t]{$n\!\! -\!\! 1$}}
\put(120,-5){\makebox(0,0)[t]{$n$}}
\put(25,20){\makebox(0,0)[l]{$0$}}
\put(120,13){\makebox(0,0)[t]{$2$}}
\end{picture}
\\
&
\\
\begin{minipage}[b]{4em}
\begin{flushright}
$C_n^{(1)}$:\\$(n \ge 2)$
\end{flushright}
\end{minipage}&
\begin{picture}(126,20)(-5,-5)
\multiput( 0,0)(20,0){3}{\circle{6}}
\multiput(100,0)(20,0){2}{\circle{6}}
\multiput(23,0)(20,0){2}{\line(1,0){14}}
\put(83,0){\line(1,0){14}}
\multiput( 2.85,-1)(0,2){2}{\line(1,0){14.3}} 
\multiput(102.85,-1)(0,2){2}{\line(1,0){14.3}} 
\multiput(59,0)(4,0){6}{\line(1,0){2}} 
\put(10,0){\makebox(0,0){$>$}}
\put(110,0){\makebox(0,0){$<$}}
\put(0,-5){\makebox(0,0)[t]{$0$}}
\put(20,-5){\makebox(0,0)[t]{$1$}}
\put(40,-5){\makebox(0,0)[t]{$2$}}
\put(100,-5){\makebox(0,0)[t]{$n\!\! -\!\! 1$}}
\put(120,-5){\makebox(0,0)[t]{$n$}}

\put(20,13){\makebox(0,0)[t]{$2$}}
\put(40,13){\makebox(0,0)[t]{$2$}}
\put(100,13){\makebox(0,0)[t]{$2$}}
\end{picture}
\\
&
\\
\begin{minipage}[b]{4em}
\begin{flushright}
$D_n^{(1)}$:\\$(n \ge 4)$
\end{flushright}
\end{minipage}&
\begin{picture}(106,40)(-5,-5)
\multiput( 0,0)(20,0){2}{\circle{6}}
\multiput(80,0)(20,0){2}{\circle{6}}
\multiput(20,20)(60,0){2}{\circle{6}}
\multiput( 3,0)(20,0){2}{\line(1,0){14}}
\multiput(63,0)(20,0){2}{\line(1,0){14}}
\multiput(39,0)(4,0){6}{\line(1,0){2}}
\multiput(20,3)(60,0){2}{\line(0,1){14}}
\put(0,-5){\makebox(0,0)[t]{$1$}}
\put(20,-5){\makebox(0,0)[t]{$2$}}
\put(80,-5){\makebox(0,0)[t]{$n\!\! - \!\! 2$}}
\put(103,-5){\makebox(0,0)[t]{$n\!\! -\!\! 1$}}
\put(25,20){\makebox(0,0)[l]{$0$}}
\put(85,20){\makebox(0,0)[l]{$n$}}
\end{picture}
\\
&
\\
$E_6^{(1)}$:&
\begin{picture}(86,60)(-5,-5)
\multiput(0,0)(20,0){5}{\circle{6}}
\multiput(40,20)(0,20){2}{\circle{6}}
\multiput(3,0)(20,0){4}{\line(1,0){14}}
\multiput(40, 3)(0,20){2}{\line(0,1){14}}
\put( 0,-5){\makebox(0,0)[t]{$1$}}
\put(20,-5){\makebox(0,0)[t]{$2$}}
\put(40,-5){\makebox(0,0)[t]{$3$}}
\put(60,-5){\makebox(0,0)[t]{$4$}}
\put(80,-5){\makebox(0,0)[t]{$5$}}
\put(45,20){\makebox(0,0)[l]{$6$}}
\put(45,40){\makebox(0,0)[l]{$0$}}
\end{picture}
\\
&
\\
$E_7^{(1)}$:&
\begin{picture}(126,40)(-5,-5)
\multiput(0,0)(20,0){7}{\circle{6}}
\put(60,20){\circle{6}}
\multiput(3,0)(20,0){6}{\line(1,0){14}}
\put(60, 3){\line(0,1){14}}
\put( 0,-5){\makebox(0,0)[t]{$0$}}
\put(20,-5){\makebox(0,0)[t]{$1$}}
\put(40,-5){\makebox(0,0)[t]{$2$}}
\put(60,-5){\makebox(0,0)[t]{$3$}}
\put(80,-5){\makebox(0,0)[t]{$4$}}
\put(100,-5){\makebox(0,0)[t]{$5$}}
\put(120,-5){\makebox(0,0)[t]{$6$}}
\put(65,20){\makebox(0,0)[l]{$7$}}
\end{picture}
\\
&
\\
$E_8^{(1)}$:&
\begin{picture}(146,40)(-5,-5)
\multiput(0,0)(20,0){8}{\circle{6}}
\put(100,20){\circle{6}}
\multiput(3,0)(20,0){7}{\line(1,0){14}}
\put(100, 3){\line(0,1){14}}
\put( 0,-5){\makebox(0,0)[t]{$0$}}
\put(20,-5){\makebox(0,0)[t]{$1$}}
\put(40,-5){\makebox(0,0)[t]{$2$}}
\put(60,-5){\makebox(0,0)[t]{$3$}}
\put(80,-5){\makebox(0,0)[t]{$4$}}
\put(100,-5){\makebox(0,0)[t]{$5$}}
\put(120,-5){\makebox(0,0)[t]{$6$}}
\put(140,-5){\makebox(0,0)[t]{$7$}}
\put(105,20){\makebox(0,0)[l]{$8$}}
\end{picture}
\\
&
\\
\end{tabular}
\begin{tabular}[t]{rl}
$F_4^{(1)}$:&
\begin{picture}(86,20)(-5,-5)
\multiput( 0,0)(20,0){5}{\circle{6}}
\multiput( 3,0)(20,0){2}{\line(1,0){14}}
\multiput(42.85,-1)(0,2){2}{\line(1,0){14.3}} 
\put(63,0){\line(1,0){14}}
\put(50,0){\makebox(0,0){$>$}}
\put( 0,-5){\makebox(0,0)[t]{$0$}}
\put(20,-5){\makebox(0,0)[t]{$1$}}
\put(40,-5){\makebox(0,0)[t]{$2$}}
\put(60,-5){\makebox(0,0)[t]{$3$}}
\put(80,-5){\makebox(0,0)[t]{$4$}}
\put(60,13){\makebox(0,0)[t]{$2$}}
\put(80,13){\makebox(0,0)[t]{$2$}}
\end{picture}
\\
&
\\
$G_2^{(1)}$:&
\begin{picture}(46,20)(-5,-5)
\multiput( 0,0)(20,0){3}{\circle{6}}
\multiput( 3,0)(20,0){2}{\line(1,0){14}}
\multiput(22.68,-1.5)(0,3){2}{\line(1,0){14.68}}
\put( 0,-5){\makebox(0,0)[t]{$0$}}
\put(20,-5){\makebox(0,0)[t]{$1$}}
\put(40,-5){\makebox(0,0)[t]{$2$}}
\put(30,0){\makebox(0,0){$>$}}
\put(40,13){\makebox(0,0)[t]{$3$}}
\end{picture}
\\
&
\\
$A^{(2)}_2$:&
\begin{picture}(26,20)(-5,-5)
\multiput( 0,0)(20,0){2}{\circle{6}}
\multiput(2.958,-0.5)(0,1){2}{\line(1,0){14.084}}
\multiput(2.598,-1.5)(0,3){2}{\line(1,0){14.804}}
\put(0,-5){\makebox(0,0)[t]{$0$}}
\put(20,-5){\makebox(0,0)[t]{$1$}}
\put(10,0){\makebox(0,0){$<$}}
\put(0,13){\makebox(0,0)[t]{$2$}}
\put(20,13){\makebox(0,0)[t]{$2$}}
\end{picture}
\\
&
\\
\begin{minipage}[b]{4em}
\begin{flushright}
$A_{2n}^{(2)}$:\\$(n \ge 2)$
\end{flushright}
\end{minipage}&
\begin{picture}(126,20)(-5,-5)
\multiput( 0,0)(20,0){3}{\circle{6}}
\multiput(100,0)(20,0){2}{\circle{6}}
\multiput(23,0)(20,0){2}{\line(1,0){14}}
\put(83,0){\line(1,0){14}}
\multiput( 2.85,-1)(0,2){2}{\line(1,0){14.3}} 
\multiput(102.85,-1)(0,2){2}{\line(1,0){14.3}} 
\multiput(59,0)(4,0){6}{\line(1,0){2}} 
\put(10,0){\makebox(0,0){$<$}}
\put(110,0){\makebox(0,0){$<$}}
\put(0,-5){\makebox(0,0)[t]{$0$}}
\put(20,-5){\makebox(0,0)[t]{$1$}}
\put(40,-5){\makebox(0,0)[t]{$2$}}
\put(100,-5){\makebox(0,0)[t]{$n\!\! -\!\! 1$}}
\put(120,-5){\makebox(0,0)[t]{$n$}}
\put(0,13){\makebox(0,0)[t]{$2$}}
\put(20,13){\makebox(0,0)[t]{$2$}}
\put(40,13){\makebox(0,0)[t]{$2$}}
\put(100,13){\makebox(0,0)[t]{$2$}}
\put(120,13){\makebox(0,0)[t]{$2$}}
\put(120,13){\makebox(0,0)[t]{$2$}}
\end{picture}
\\
&
\\
$A^{(2)\dagger}_2$:&
\begin{picture}(26,20)(-5,-5)
\multiput( 0,0)(20,0){2}{\circle{6}}
\multiput(2.958,-0.5)(0,1){2}{\line(1,0){14.084}}
\multiput(2.598,-1.5)(0,3){2}{\line(1,0){14.804}}
\put(0,-5){\makebox(0,0)[t]{$0$}}
\put(20,-5){\makebox(0,0)[t]{$1$}} \put(10,0){\makebox(0,0){$>$}}
\put(0,13){\makebox(0,0)[t]{$$}} \put(20,13){\makebox(0,0)[t]{$$}}
\end{picture}
\\
&
\\
\begin{minipage}[b]{4em}
\begin{flushright}
$A_{2n}^{(2)\dagger}$:\\$(n \ge 2)$
\end{flushright}
\end{minipage}&
\begin{picture}(126,20)(-5,-5)
\multiput( 0,0)(20,0){3}{\circle{6}}
\multiput(100,0)(20,0){2}{\circle{6}}
\multiput(23,0)(20,0){2}{\line(1,0){14}}
\put(83,0){\line(1,0){14}}
\multiput( 2.85,-1)(0,2){2}{\line(1,0){14.3}} 
\multiput(102.85,-1)(0,2){2}{\line(1,0){14.3}} 
\multiput(59,0)(4,0){6}{\line(1,0){2}} 
\put(10,0){\makebox(0,0){$>$}} \put(110,0){\makebox(0,0){$>$}}
\put(0,-5){\makebox(0,0)[t]{$0$}}
\put(20,-5){\makebox(0,0)[t]{$1$}}
\put(40,-5){\makebox(0,0)[t]{$2$}}
\put(100,-5){\makebox(0,0)[t]{$n\!\! -\!\! 1$}}
\put(120,-5){\makebox(0,0)[t]{$n$}}
\put(0,13){\makebox(0,0)[t]{$$}}
\put(20,13){\makebox(0,0)[t]{$$}}
\put(40,13){\makebox(0,0)[t]{$$}}
\put(100,13){\makebox(0,0)[t]{$$}}
\put(120,13){\makebox(0,0)[t]{$$}}
\put(120,13){\makebox(0,0)[t]{$$}}
\end{picture}
\\
&
\\
\begin{minipage}[b]{4em}
\begin{flushright}
$A_{2n-1}^{(2)}$:\\$(n \ge 3)$
\end{flushright}
\end{minipage}&
\begin{picture}(126,40)(-5,-5)
\multiput( 0,0)(20,0){3}{\circle{6}}
\multiput(100,0)(20,0){2}{\circle{6}}
\put(20,20){\circle{6}}
\multiput( 3,0)(20,0){3}{\line(1,0){14}}
\multiput(83,0)(20,0){1}{\line(1,0){14}}
\put(20,3){\line(0,1){14}}
\multiput(102.85,-1)(0,2){2}{\line(1,0){14.3}} 
\multiput(59,0)(4,0){6}{\line(1,0){2}} 
\put(110,0){\makebox(0,0){$<$}}
\put(0,-5){\makebox(0,0)[t]{$1$}}
\put(20,-5){\makebox(0,0)[t]{$2$}}
\put(40,-5){\makebox(0,0)[t]{$3$}}
\put(100,-5){\makebox(0,0)[t]{$n\!\! -\!\! 1$}}
\put(120,-5){\makebox(0,0)[t]{$n$}}
\put(25,20){\makebox(0,0)[l]{$0$}}

\put(120,13){\makebox(0,0)[t]{$2$}}
\end{picture}
\\
&
\\
\begin{minipage}[b]{4em}
\begin{flushright}
$D_{n+1}^{(2)}$:\\$(n \ge 2)$
\end{flushright}
\end{minipage}&
\begin{picture}(126,20)(-5,-5)
\multiput( 0,0)(20,0){3}{\circle{6}}
\multiput(100,0)(20,0){2}{\circle{6}}
\multiput(23,0)(20,0){2}{\line(1,0){14}}
\put(83,0){\line(1,0){14}}
\multiput( 2.85,-1)(0,2){2}{\line(1,0){14.3}} 
\multiput(102.85,-1)(0,2){2}{\line(1,0){14.3}} 
\multiput(59,0)(4,0){6}{\line(1,0){2}} 
\put(10,0){\makebox(0,0){$<$}}
\put(110,0){\makebox(0,0){$>$}}
\put(0,-5){\makebox(0,0)[t]{$0$}}
\put(20,-5){\makebox(0,0)[t]{$1$}}
\put(40,-5){\makebox(0,0)[t]{$2$}}
\put(100,-5){\makebox(0,0)[t]{$n\!\! -\!\! 1$}}
\put(120,-5){\makebox(0,0)[t]{$n$}}

\put(20,13){\makebox(0,0)[t]{$2$}}
\put(40,13){\makebox(0,0)[t]{$2$}}
\put(100,13){\makebox(0,0)[t]{$2$}}
\end{picture}
\\
&
\\
$E_6^{(2)}$:&
\begin{picture}(86,20)(-5,-5)
\multiput( 0,0)(20,0){5}{\circle{6}}
\multiput( 3,0)(20,0){2}{\line(1,0){14}}
\multiput(42.85,-1)(0,2){2}{\line(1,0){14.3}} 
\put(63,0){\line(1,0){14}}
\put(50,0){\makebox(0,0){$<$}}
\put( 0,-5){\makebox(0,0)[t]{$0$}}
\put(20,-5){\makebox(0,0)[t]{$1$}}
\put(40,-5){\makebox(0,0)[t]{$2$}}
\put(60,-5){\makebox(0,0)[t]{$3$}}
\put(80,-5){\makebox(0,0)[t]{$4$}}

\put(60,13){\makebox(0,0)[t]{$2$}}
\put(80,13){\makebox(0,0)[t]{$2$}}
\end{picture}
\\
&
\\
$D_4^{(3)}$:&
\begin{picture}(46,20)(-5,-5)
\multiput( 0,0)(20,0){3}{\circle{6}}
\multiput( 3,0)(20,0){2}{\line(1,0){14}}
\multiput(22.68,-1.5)(0,3){2}{\line(1,0){14.68}}
\put( 0,-5){\makebox(0,0)[t]{$0$}}
\put(20,-5){\makebox(0,0)[t]{$1$}}
\put(40,-5){\makebox(0,0)[t]{$2$}}
\put(30,0){\makebox(0,0){$<$}}

\put(40,13){\makebox(0,0)[t]{$3$}}
\end{picture}
\\
&
\\
\end{tabular}
\caption{Dynkin diagrams for $X^{(r)}_N$.
The enumeration of the nodes with
$I = \{0,1,\ldots, n\}$ is specified under or the right side of the nodes.
In addition, the numbers $t_i$ (resp. $t^\vee_i$) defined in
\eqref{eq:ttdef} are attached {\em above} the nodes for $r=1$ (resp. $r>1$)
if and only if $t_i \neq 1$ (resp. $t^\vee_i \neq 1$).\label{fig:Dynkin}}
\end{figure}
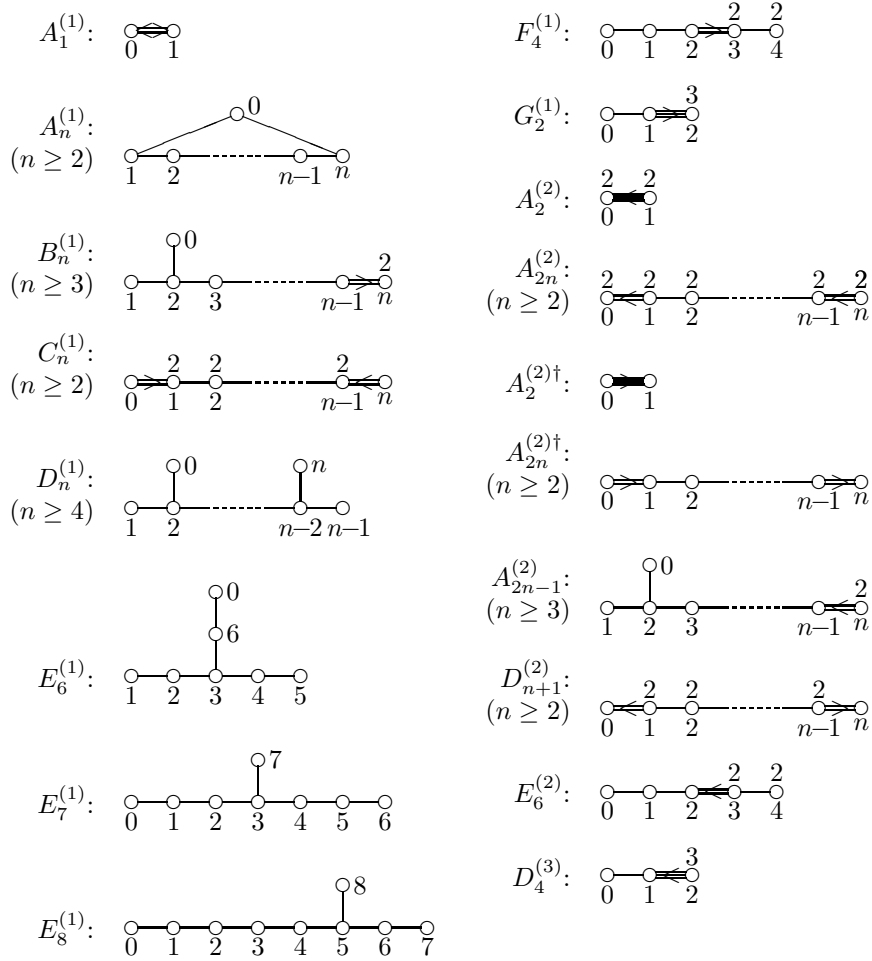}

Every affine algebra $\geh$ has a simple Lie subalgebra
$\gehb$ obtained by removing the 0-node from the Dynkin diagram. This is
summarized in the following table:
\begin{equation} \label{eq:subalgs}
\begin{array}{c|ccccccc}
\geh & X_n^{(1)} & A_{2n}^{(2)} & A_{2n}^{(2)\dagger} & A_{2n-1}^{(2)} &
D_{n+1}^{(2)} & E_6^{(2)} & D_4^{(3)} \\ \hline
\gehb & X_n & C_n & B_n & C_n & B_n & F_4 & G_2
\end{array}
\end{equation}

Let $\alpha_i,h_i,\La_i$ ($i \in \I$) be the simple roots, simple
coroots, and fundamental weights of $\geh$. Let $\delta$ and $c$
denote the generator of imaginary roots and the canonical central
element, respectively. Recall that $\delta=\sum_{i\in \I}a_i\alpha_i$
and $c = \sum_{i \in \I}a^\vee_i h_i$, where the Kac
labels $a_i$ are the unique set of relatively prime positive
integers giving the linear dependency of the columns of the Cartan
matrix $A$ (that is, $A (a_0,\dotsc,a_n)^t = 0$). Explicitly,
\begin{equation}
\delta =
\begin{cases}
\alpha_0+\cdots + \alpha_n & \text{if $\geh=A^{(1)}_n$} \\
\alpha_0+\alpha_1+2\alpha_2+\dotsm+2\alpha_n & \text{if $\geh=B^{(1)}_{n}$} \\
\alpha_0+2\alpha_1+\dotsm+2\alpha_{n-1}+\alpha_n & \text{if $\geh=C^{(1)}_{n}$} \\
\alpha_0+\alpha_1+2\alpha_2+\dotsm+2\alpha_{n-2}
+\alpha_{n-1}+\alpha_n
& \text{if $\geh=D^{(1)}_{n}$} \\
\alpha_0+\alpha_1+2\alpha_2+3\alpha_3+2\alpha_4+\alpha_5+2\alpha_6
& \text{if $\geh=E_6^{(1)}$} \\
\alpha_0+2\alpha_1+3\alpha_2+4\alpha_3+3\alpha_4+2\alpha_5+\alpha_6+\alpha_7
& \text{if $\geh=E_7^{(1)}$} \\
\alpha_0+2\alpha_1+3\alpha_2+4\alpha_3+5\alpha_4+6\alpha_5+4\alpha_6+2\alpha_7+3\alpha_8
& \text{if $\geh=E_8^{(1)}$} \\
\alpha_0+2\alpha_1+3\alpha_2+4\alpha_3+2\alpha_4
& \text{if $\geh=F_4^{(1)}$} \\
\alpha_0+2\alpha_1+3\alpha_2
& \text{if $\geh=G_2^{(1)}$} \\
2\alpha_0+2\alpha_1+\dotsm+2\alpha_{n-1}+\alpha_n
& \text{if $\geh=A^{(2)}_{2n}$} \\
\alpha_0+2\alpha_1+\dotsm+2\alpha_{n-1}+2\alpha_n
& \text{if $\geh=A^{(2)\dagger}_{2n}$} \\
\alpha_0+\alpha_1+2\alpha_2+\dotsm+2\alpha_{n-1}+\alpha_n
& \text{if $\geh=A^{(2)}_{2n-1}$} \\
\alpha_0+\alpha_1+\dotsm+\alpha_{n-1}+\alpha_n
&\text{if $\geh=D^{(2)}_{n+1}$} \\
\alpha_0+2\alpha_1+3\alpha_2+2\alpha_3+\alpha_4
& \text{if $\geh=E_6^{(2)}$} \\
\alpha_0+2\alpha_1+\alpha_2
& \text{if $\geh=D_4^{(3)}$}.
\end{cases}
\end{equation}
The dual Kac label $a^\vee_i$ is the label $a_i$ for the affine
Dynkin diagram obtained by ``reversing the arrows" of the Dynkin
diagram of $\geh$, or equivalently, the coefficients giving the
linear dependency of the rows of the Cartan matrix $A$.

Let $P=\bigoplus_{a\in \I} \Z \La_a\oplus \Z \delta$ be the
weight lattice of $\geh$ and $P^+=\bigoplus_{a\in\I} \Z_{\ge 0} \La_a$.
Similarly, let $\Pfin=\bigoplus_{a\in\J} \Z\Lab_a$ be the weight lattice
of $\gehb$, $\Pfin^+=\bigoplus_{a\in\J} \Z_{\ge 0}\Lab_a$,
$\Qfin=\bigoplus_{a\in\J} \Z \alpha_a$ the root lattice of $\gehb$
and $\Qfin^+=\bigoplus_{a\in\J} \Z_{\ge 0} \alpha_a$
with simple roots and fundamental weights $\alpha_a,\Lab_a$ for $a\in\J$.
For $\la,\mu\in \Pfin$ write $\la\dom\mu$ if $\la-\mu\in \Qfin^+$.

For $i \in \I$ let
\begin{equation}\label{eq:ttdef}
t_i = \max(\frac{a_i}{a^\vee_i},a^\vee_0), \qquad
t^\vee_i = \max(\frac{a^\vee_i}{a_i},a_0).
\end{equation}
The values $t_i$ are given in Figure \ref{fig:Dynkin}. We shall
only use $t^\vee_i$ and $t_i$ for $i \in \J$. For $a\in\J$ we have
\begin{equation*}
t^\vee_a = 1 \,\,\text{ if $r = 1$,} \qquad
t_a = a_0^\vee \,\,\text{ if $r > 1$}.
\end{equation*}

Let $(\cdot|\cdot)$ be the normalized invariant form on $P$
\cite{Kac}. It satisfies
\begin{equation}
  (\alpha_i|\alpha_j) = \dfrac{a_i^\vee}{a_i} A_{ij}
\end{equation}
for $i,j\in I$. In particular
\begin{equation} \label{eq:form norm}
  (\alpha_a|\alpha_a)=\dfrac{2r}{a_0^\vee}
\end{equation}
if $\alpha_a$ is a long root.

\subsection{Crystals}
The quantized universal enveloping algebra $U_q(\geh)$ associated
with a symmetrizable Kac--Moody Lie algebra $\geh$ was introduced
independently by Drinfeld \cite{D:1985} and Jimbo \cite{J:1985} in their
study of two dimensional solvable lattice models in statistical
mechanics. The parameter $q$ corresponds to the temperature of the
underlying model. Kashiwara \cite{K:1990} showed that at zero
temperature or $q=0$ the representations of $U_q(\geh)$ have
bases, which he coined crystal bases, with a beautiful
combinatorial structure and favorable properties such as
uniqueness and stability under tensor products.

Let $\geh'$ be the derived subalgebra of $\geh$. Denote the
corresponding quantized universal enveloping algebras of
$\geh\supset \geh'\supset \gehb$ by $U_q(\geh)\supset
U'_q(\geh)\supset U_q(\gehb)$.

In \cite{HKOTY:1999,HKOTT:2002} it is conjectured that there is a family of
finite-dimensional irreducible $U'_q(\geh)$-modules
$\{W^{(a)}_i\mid a\in\J,i\in \Z_{>0}\}$ which, unlike most
finite-dimensional $U'_q(\geh)$-modules, have crystal bases
$B^{a,i}$. This family is conjecturally characterized in several
different ways:
\begin{enumerate}
\item Its characters form the unique solutions of a system
of quadratic relations (the $Q$-system) \cite{KR:1987}.
\item Every crystal graph of an irreducible integrable finite-dimensional
$U'_q(\geh)$-module, is a tensor product of the $B^{a,i}$.
\item For $\la\in P$ let $V(\la)$ be the universal extremal weight module
defined in \cite[Section 3]{K:2000} and $B(\la)$ its crystal base,
with unique vector $u_\la\in B(\la)$ of weight $\la$. Then the
affinization of $B^{a,i}$ (in the sense of \cite{KKMMNN:1992a}) is
isomorphic to the connected component of $u_\la$ in $B(\la)$, for
the weight $\la=i\Lab_a$.
\end{enumerate}

In light of point (2) above, we consider the category of crystal
graphs given by tensor products of the crystals $B^{a,i}$.

We introduce notation for tensor products of $B^{a,i}$. Let
\begin{equation}\label{eq:B}
 B = \bigotimes_{(a,i)\in\J\times\Z_{>0}}
(B^{a,i})^{\otimes L_i^{(a)}},
\end{equation}
where only finitely many $L_i^{(a)}$ are nonzero.
In type $A^{(1)}_n$ this is the tensor product of modules, which,
when restricted to $A_n$, are irreducible modules indexed by
rectangular partitions. The set of classically restricted paths
(or classical highest weight vectors) in $B$ of weight
$\la\in\Pfin^+$ is by
definition
\begin{equation*}
\Path(B,\la)=\{b\in B\mid \text{$\wt(b)=\la$ and $e_ib$
undefined for all $i\in \J$} \}.
\end{equation*}
Here $e_i$ is given by the crystal graph. For $b,b'\in B^{a,i}$
we have $b'=e_i(b)$ if there is an arrow $b'\stackrel{i}{\longrightarrow}b$
in the crystal graph; if no such arrow exists then $e_i(b)$ is undefined.
Similarly, $b'=f_i(b)$ if there is an arrow $b\stackrel{i}{\longrightarrow}b'$
in the crystal graph; if no such arrow exists then $f_i(b)$ is undefined.
If $B_1$ and $B_2$ are crystals, then for $b_1\otimes b_2\in B_1\otimes B_2$
the action of $e_i$ is defined as
\begin{equation*}
e_i(b_1\otimes b_2)=\begin{cases}
e_ib_1 \otimes b_2 &\text{if $\varepsilon_i(b_1)>\varphi_i(b_2)$,}\\
b_1\otimes e_i b_2 &\text{else,}
\end{cases}
\end{equation*}
where $\varepsilon_i(b)=\max\{k\mid e_i^k b \;\text{is defined}\}$ and
$\varphi_i(b)=\max\{k\mid f_i^k b\;\text{is defined}\}$. This is the opposite
of the notation used by Kashiwara \cite{K:1990}.

\subsection{Simple crystals}
Let $W$ be the Weyl group of $\geh$, $\{s_i\mid i\in \I\}$ the simple
reflections in $W$. Let $B$ be the crystal graph of an integrable
$U_q(\geh)$-module. Say that $b\in B$ is an extremal vector of weight
$\la\in P$ provided that $\wt(b)=\la$ and there exists a family of
elements $\{b_w\mid w\in W \}\subset B$ such that
\begin{enumerate}
\item $b_w=b$ for $w=e$.
\item If $\inner{h_i}{w\la}\ge 0$ then $e_i(b_w)=\es$ and
$f_i^{\inner{h_i}{w\la}}(b_w)=b_{s_i w}$.
\item If $\inner{h_i}{w\la}\le 0$ then $f_i(b_w)=\es$ and
$e_i^{\inner{h_i}{w\la}}(b_w)=b_{s_i w}$.
\end{enumerate}

Following \cite{AK:1997}, say that a
$U'_q(\geh)$-crystal $B$ is \textit{simple} if
\begin{enumerate}
\item $B$ is the crystal base of a finite dimensional integrable
$U'_q(\geh)$-module.
\item There is a weight $\la\in \Pfin^+$ such that $B$ has a
unique vector (denoted $u(B)$) of weight $\la$, and the weight of
any extremal vector of $B$ is contained in $\Wb\la$ where $\Wb$ is the
Weyl group of $\gehb$.
\end{enumerate}

In the definition of simple crystal in \cite{AK:1997}, condition 1 is not
present. However we always want to assume both conditions,
so it is convenient to include condition 1 in the definition above.

\begin{theorem}[\cite{AK:1997}]\mbox{}
\begin{enumerate}
\item Simple crystals are connected.
\item The tensor product of simple crystals is simple.
\end{enumerate}
\end{theorem}

For the $U'_q(\geh)$-crystal $B$,
define $\epsilon,\varphi:B\rightarrow P$ by
\begin{equation*}
\epsilon(b)=\sum_{i\in \I} \epsilon_i(b) \La_i \qquad \text{and}\qquad
\varphi(b)=\sum_{i\in \I} \varphi_i(b) \La_i.
\end{equation*}
Then the level of $B$ is
\begin{equation}
\lev(B)=\min \{\inner{c}{\epsilon(b)}\mid b\in B\}.
\end{equation}

\subsection{Dual crystals}
The notion of a dual crystal is given in
\cite[Section 7.4]{K:1995}. Let $B$ be a $U_q(\geh)$-crystal. Then
there is a $U_q(\geh)$-crystal denoted $B^\vee$ obtained from $B$
by reversing arrows. That is, $B^\vee=\{b^\vee\mid b\in B\}$ with
\begin{equation} \label{eq:dual crystal}
\begin{split}
\wt(b^\vee)&=-\wt(b) \\
\epsilon_i(b^\vee)&=\varphi_i(b) \\
\varphi_i(b^\vee)&=\epsilon_i(b) \\
e_i(b^\vee) &= (f_i(b))^\vee \\
f_i(b^\vee) &= (e_i(b))^\vee.
\end{split}
\end{equation}

\begin{prop} \label{pp:dual tensor} \cite{K:1995} There is an isomorphism
$(B_2\otimes B_1)^\vee \cong B_1^\vee \otimes B_2^\vee$ given by
$(b_2\otimes b_1)^\vee \mapsto b_1^\vee \otimes b_2^\vee$.
\end{prop}

\subsection{One dimensional sums}
\label{subsec:1dsum}
In this section we recall the structure of a $U'_q(\geh)$-crystal
as a graded $U_q(\gehb)$-crystal. The grading is given by the intrinsic
energy function $D:B\rightarrow\Z$. For $b\in B$, one may define
$D(b)$ as the minimum number of times $e_0$ occurs in a sequence
of operators involving $e_i,f_i$ for $i\in\J$ and $e_0$,
leading from $u(B)$ to $b$.
However we prefer to work with the following concrete definition
when $B$ is a tensor product of crystals of the form $B^{r,s}$.
This definition essentially comes from \cite{HKOTT:2002},
but it is useful to formulate it as follows \cite{OSS:2001}.

Let $B_1, B_2$ be simple $U'_q(\geh)$-crystals. It was shown in
\cite[Section 4]{KKMMNN:1992a} that there is a unique isomorphism of
$U'_q(\geh)$-crystals
$\R=\R_{B_2,B_1}:B_2\otimes B_1\rightarrow B_1\otimes B_2$, called
the combinatorial $\R$ matrix. In addition there exists a
function $H:B_1\otimes B_2\to \Z$ called the local energy function,
that is unique up to a global additive constant,
which is constant on $\J$ components and satisfies for all $b_2\in B_2$
and $b_1\in B_1$ with $\R(b_2\otimes b_1)=b_1'\otimes b_2'$
\begin{equation} \label{eq:H}
  H(e_0(b_2\otimes b_1))=
  H(b_2\otimes b_1)+
  \begin{cases}
    -1 & \text{if $\epsilon_0(b_2)>\varphi_0(b_1)$ and
    $\epsilon_0(b_1')>\varphi_0(b_2')$} \\
    1 & \text{if $\epsilon_0(b_2)\le\varphi_0(b_1)$ and
    $\epsilon_0(b_1')\le \varphi_0(b_2')$} \\
    0 & \text{otherwise.}
  \end{cases}
\end{equation}
We shall normalize the local energy function by the condition
$H(u(B_2)\otimes u(B_1))=0$.

It was conjectured in \cite{HKOTT:2002} that
\begin{equation} \label{eq:nat}
\varphi(b^\natural)=\lev(B^{r,s}) \La_0
\qquad\text{for a unique $b^\natural\in B^{r,s}$.}
\end{equation}
For a given crystal $B^{r,s}$, denote this element also by
$u^\natural(B^{r,s})$.
Define the function $D_{B^{r,s}}:B^{r,s}\rightarrow\Z$ by
\begin{equation} \label{eq:DBrs}
  D_{B^{r,s}}(b) = H(b\otimes b^\natural)-H(u(B^{r,s})\otimes b^\natural)
\end{equation}
where $H=H_{B^{r,s},B^{r,s}}$ is the local energy function.
In all cases in which the $U'_q(\geh)$-module $W^{(r)}_s$
and its crystal base $B^{r,s}$ have been
constructed, \eqref{eq:nat} holds and \eqref{eq:DBrs}
agrees with the explicit grading on $B^{r,s}$ specified
in a case-by-case manner in the appendices
of~\cite{HKOTY:1999,HKOTT:2002}.

A graded simple crystal $(B,D)$ is a simple crystal $B$ together
with a function $D:B\rightarrow\Z$.
Let $(B_j,D_j)$ be a graded simple $U'_q(\geh)$-crystal
and $u_j=u(B_j)$, for $1\le j\le L$.
Let $B=B_L\otimes\dotsm\otimes B_1$.
Following \cite{NY:1997} define the energy function
$E_B:B\rightarrow\Z$ by
\begin{equation}\label{eq:NY}
E_B = \sum_{1\le i<j\le L} H_i\R_{i+1}\R_{i+2}\dotsm\R_{j-1},
\end{equation}
where $\R_i$ is the combinatorial $\R$-matrix and $H_i$ is the local
energy function, where the subscript $i$ indicates that the operators
act on the $i$-th and $(i+1)$-st tensor factors from the right.
This given, define $D_B,D'_B:B\rightarrow\Z$ by
\begin{equation}\label{eq:D}
\begin{split}
  D'_B &= E_B + \sum_{j=1}^L D_j \R_1\R_2\dotsm\R_{j-1} \\
  D_B(b) &= D'_B(b) - D'_B(u(B))
\end{split}
\end{equation}
where $D_j:B_j\rightarrow\Z$ acts on the rightmost tensor factor.
Then we say that the graded simple crystal $(B,D_B)$ is the
tensor product of the graded simple crystals $(B_j,D_j)$.

\begin{theorem} \cite{OSS:2001} Graded simple $U'_q(\geh)$-crystals
form a tensor category.
\end{theorem}

Now suppose that for all $j$, $B_j$ has the form $B^{r,s}$
and $D_j=D_{B^{r,s}}$ is the intrinsic energy as defined above.
Then the function $D_B$ is called the intrinsic energy of $B$.
Let $b_j^\natural\in B_j$ be as in \eqref{eq:nat}. Then conjecturally
there is an element $b^\natural\in B$ such that $b_j^\natural$
is the leftmost tensor factor in $R_{L-1} \dotsm R_{j+1} R_j b^\natural$.

Using the Yang-Baxter equation for $R$ and the fact that
$R_{B\otimes B}$ is the identity for any $B$, it follows that \cite{HKOTT:2002}
\begin{equation} \label{eq:b nat}
  D_B(b) = H(b\otimes b^\natural)-H(u(B)\otimes b^\natural).
\end{equation}

The one-dimensional sum $X(B,\la;q)\in\Z[q,q^{-1}]$ is the generating
function of paths graded by the intrinsic energy
\begin{equation} \label{eq:onedim}
  X(B,\la;q) = \sum_{b\in \Path(B,\la)} q^{D_B(b)}.
\end{equation}

\subsection{Crystals of type $B_n,C_n,D_n$}
\label{subsec:clcrystals}

In this and the next section we will describe the
classical highest weight crystals $B(s\Lab_1)$
and the finite dimensional affine crystals $B^{1,s}$ for
all nonexceptional types
as weakly increasing words $b$ in an alphabet $\X$.
They are also determined by
$x(b)=(x_i)_{i\in \X}$ where $x_i$ is the number of $i$'s
in $b$. Whenever an operation yields a negative value for an
$x_i$ it will be undefined.

According to \cite{KN:1994}, the crystal of $B(\Lab_1)$ has underlying set
\begin{align*}
  \X &= \{1<2<\cdots<n<\0<\bar{n}<\cdots<\bar{2}<\bar{1} \}
     && \text{for $B_n$} \\
  \X &= \{1<2<\cdots<n<\bar{n}<\cdots<\bar{2}<\bar{1} \}
     && \text{for $C_n$} \\
  \X &= \{1<2<\cdots<\begin{matrix}{n} \\ {\bar{n}} \end{matrix} <\cdots<\bar{2}<\bar{1} \}
     && \text{for $D_n$.}
\end{align*}
The crystal $B(s\Lab_1)$ is the set of weakly increasing words of length
$s$ in the alphabet $\X$ such that, in addition, for type $B_n$ there is
at most one $\0$, and in type $D_n$, there are either no letters $n$ or
no letters $\bar{n}$.

The crystal operators $e_i$ on $B(s\Lab_1)$ are given by
\begin{equation} \label{eq:BCe}
e_i b = \begin{cases}
 (x_1,\ldots,x_i+1,x_{i+1}-1,\ldots,\xb_1) & \text{if $x_{i+1}>\xb_{i+1}$}\\
 (x_1,\ldots,\xb_{i+1}+1,\xb_i-1,\ldots,\xb_1) & \text{if $x_{i+1}\le\xb_{i+1}$}
\end{cases}
\end{equation}
with the following exceptions:
\begin{equation}
\begin{aligned}
\text{Type $B_n$:}& & e_n b &= \begin{cases}
 (x_1,\ldots,x_n,x_\0+1,\xb_n-1,\ldots,\xb_1) & \text{if $x_\0=0$}\\
 (x_1,\ldots,x_n+1,x_\0-1,\xb_n,\ldots,\xb_1) & \text{if $x_\0=1$.}
\end{cases}  \\
\text{Type $C_n$:}& & e_n b &= (x_1,\ldots,x_n+1,\xb_n-1,\ldots,\xb_1) \\
\text{Type $D_n$:\;}& & e_{n-1} b &= \begin{cases}
 (x_1,\ldots,x_{n-1}+1,x_n-1,\xb_n,\ldots,\xb_1) & \text{if $x_n>0$}\\
 (x_1,\ldots,x_n,\xb_n+1,\xb_{n-1}-1,\ldots,\xb_1) & \text{if $x_n=0$}
\end{cases}\\
& & e_n b &= \begin{cases}
 (x_1,\ldots,x_{n-1}+1,x_n,\xb_n-1,\ldots,\xb_1) & \text{if $\xb_n>0$}\\
 (x_1,\ldots,x_n+1,\xb_n,\xb_{n-1}-1,\ldots,\xb_1) & \text{if $\xb_n=0$.}
\end{cases}
\end{aligned}
\end{equation}

\subsection{Affine crystals $B^{1,s}$}

We recall the crystals $B^{1,s}$ from \cite{KKMMNN:1992}
(and \cite{KKM:1994} for type $C_n^{(1)}$).
The affine algebra $\geh$ has simple Lie subalgebra of type given
in \eqref{eq:subalgs}. There is an isomorphism of classical crystals
\begin{equation}
B^{1,s} \cong
\begin{cases}
  B(s\Lab_1) & \text{for types $B_n^{(1)},D_n^{(1)},A_{2n-1}^{(2)}$} \\[2mm]
  \displaystyle{\bigoplus_{s'\le s} B(s'\Lab_1)} & \text{for types $A_{2n}^{(2)},D_{n+1}^{(2)}$} \\[2mm]
  \displaystyle{\bigoplus_{\substack{s'\le s \\ s-s'\in 2\Z}} B(s'\Lab_1)} &
  \text{for type $C_n^{(1)},A_{2n}^{(2)\dagger}$.}
\end{cases}
\end{equation}
The crystal operators $e_i$ for $1\le i\le n$ are given in
subsection \ref{subsec:clcrystals}. The operator $e_0$ is given by
\begin{equation}
\begin{aligned}[3]
&\text{Type $B_n^{(1)},D_n^{(1)},A_{2n-1}^{(2)}$:\;} & e_0 b &= \begin{cases}
 (x_1,x_2-1,\ldots,\xb_2,\xb_1+1) & \text{if $x_2>\xb_2$}\\
 (x_1-1,x_2,\ldots,\xb_2+1,\xb_1) & \text{if $x_2\le \xb_2$}
\end{cases} \\
&\text{Type $A_{2n}^{(2)},D_{n+1}^{(2)}$:} & e_0 b &= \begin{cases}
  (x_1-1,x_2,\dots,\xb_2,\xb_1) & \text{if $x_1> \xb_1$} \\
  (x_1,x_2,\dots,\xb_2,\xb_1+1) &\text{if $x_1 \le \xb_1$}
\end{cases} \\
&\text{Type $C_n^{(1)},A_{2n}^{(2)\dagger}$:} & e_0 b &= \begin{cases}
  (x_1-2,x_2,\dots,\xb_2,\xb_1) & \text{if $x_1\ge \xb_1+2$} \\
  (x_1-1,x_2,\dots,\xb_2,\xb_1+1) &\text{if $x_1=\xb_1+1$} \\
  (x_1,x_2,\dots,\xb_2,\xb_1+2) &\text{if $x_1 \le \xb_1$.}
\end{cases}
\end{aligned}
\end{equation}

\section{Virtual crystals}
\label{sec:virtual}

\subsection{Embeddings of affine algebras}

As given in \eqref{eq:embed}, there are natural inclusions of the
affine Lie algebras. These embeddings do not carry over to the
corresponding quantum algebras. Nevertheless we expect that such
embeddings exist for crystals. Note that every affine algebra can
be embedded into one of type $A^{(1)}, D^{(1)}$ and $E^{(1)}$
which are the untwisted affine algebras whose canonical simple Lie
subalgebra is simply-laced. Crystal embeddings $C_n^{(1)},
A_{2n}^{(2)}, A_{2n}^{(2)\dagger},D_{n+1}^{(2)} \hookrightarrow
A_{2n-1}^{(1)}$ are studied in \cite{OSS:2001}.

Consider one of the embeddings given in \eqref{eq:embed}
of an affine algebra with Dynkin diagram $X$ into one
with diagram $Y$. We consider a graph automorphism $\aut$ of $Y$
that fixes the 0 node. For type $A_{2n-1}^{(1)}$, $\aut(i)=2n-i$ (mod $2n$).
For type $D_{n+1}^{(1)}$ the automorphism interchanges the nodes $n$ and
$n+1$ and fixes all other nodes. There is an additional
automorphism for type $D_4^{(1)}$, namely, the cyclic permutation
of the nodes 1,2 and 3. For type $E_6^{(1)}$ the automorphism
exchanges nodes 1 and 5 and nodes 2 and 4.

Let $\I^X$ and $\I^Y$ be the vertex sets of the diagrams $X$ and $Y$
respectively, $\I^Y/\aut$ the set of orbits of the action of
$\aut$ on $\I^Y$, and $\bij:\I^X\rightarrow \I^Y/\aut$ a bijection
which preserves edges and sends $0$ to $0$.

\begin{example}
If $X$ is one of $C_n^{(1)}, A_{2n}^{(2)},
A_{2n}^{(2)\dagger},D_{n+1}^{(2)}$ and $Y=A_{2n-1}^{(1)}$, then
$\bij(0)=0$, $\bij(i)=\{i,2n-i\}$ for $1\le i<n$ and $\bij(n)=n$.

If $X=B_n^{(1)}$ or $A_{2n-1}^{(2)}$ and $Y=D_{n+1}^{(1)}$,
then $\bij(i)=i$ for $i<n$ and $\bij(n)=\{n,n+1\}$.

If $X$ is $E_6^{(2)}$ or $F_4^{(1)}$ and $Y=E_6^{(1)}$, then
$\bij(0)=0$, $\bij(1)=1$, $\bij(2)=3$, $\bij(3)=\{2,4\}$ and
$\bij(4)=\{1,5\}$.

If $X$ is $D_4^{(3)}$ or $G_2^{(1)}$ and $Y=D_4^{(1)}$, then
$\bij(0)=0$, $\bij(1)=2$ and $\bij(2)=\{1,3,4\}$.
\end{example}

To describe the embedding we endow the bijection $\bij$ with
additional data. For each $i\in I^X$ we shall define a
multiplication factor $\mult_i$ that depends on the location of
$i$ with respect to a distinguished arrow (multiple bond)
in $X$. Removing the arrow leaves two connected components. The
factor $\mult_i$ is defined as follows:
\begin{enumerate}
\item Suppose $X$ has a unique arrow.
\begin{enumerate}
\item Suppose the arrow points towards the component of $0$.
Then $\mult_i=1$ for all $i\in I^X$.
\item Suppose the arrow points away from the component of $0$.
Then $\mult_i$
is the order of $\aut$ for $i$ in the component of $0$ and is
$1$ otherwise.
\end{enumerate}
\item Suppose $X$ has two arrows, that is, $Y=A_{2n-1}^{(1)}$.
Then $\mult_i=1$ for $1\le i\le n-1$. For
$i\in\{0,n\}$, $\mult_i=2$ (which is the order of $\aut$)
if the arrow incident to $i$ points away from it
and is $1$ otherwise.
\end{enumerate}

\begin{example}
For $X=B_n^{(1)}$ and $Y=D_{n+1}^{(1)}$ we have $\mult_i=2$ if $0\le i \le n-1$
and $\mult_n=1$.
For $X=A_{2n-1}^{(2)}$ and $Y=D_n^{(1)}$ we have $\mult_i=1$ for all $i$.
\end{example}

The embedding $\emb:P^X\to P^Y$ of weight lattices is defined by
\begin{equation*}
\emb(\La^X_i) = \mult_i \sum_{j\in \bij(i)} \La^Y_j.
\end{equation*}
As a consequence we have
\begin{equation*}
\begin{split}
  \emb(\alpha^X_i) &= \mult_i \sum_{j\in \bij(i)} \alpha^Y_j \\
  \emb(\delta^X) &= a_0 \mult_0 \,\delta^Y.
\end{split}
\end{equation*}

\subsection{Virtual crystals}

Suggested by the embeddings $X\hookrightarrow Y$ of affine algebras,
we wish to realize crystals of type $X$ using crystals of type $Y$.

Let $\Vh$ be a $Y$-crystal. We define the virtual crystal
operators $\eh_i,\fh_i$ for $i\in I^X$ as the composites of $Y$-crystal operators
$f_j,e_j$ given by
\begin{equation*}
\begin{split}
\fh_i &= \prod_{j\in \bij(i)} f_j^{\mult_i}\\
\eh_i &= \prod_{j\in \bij(i)} e_j^{\mult_i}.
\end{split}
\end{equation*}
These are designed to simulate $X$-crystal operators $f_i,e_i$ for $i\in \I^X$.
The type $Y$ operators on the right hand side, may be performed
in any order, since distinct nodes $j,j'\in \bij(i)$ are
not adjacent in $Y$ and thus their corresponding raising and
lowering operators commute.

A virtual crystal is a pair $(\V,\Vh)$ such that:
\begin{enumerate}
\item $\Vh$ is a $Y$-crystal.
\item $\V\subset \Vh$ is closed under $\eh_i,\fh_i$ for $i\in I^X$.
\item There is an $X$-crystal $B$ and an $X$-crystal isomorphism
$\viso:B\rightarrow \V$ such that $e_i,f_i$ correspond to $\eh_i,\fh_i$.
\end{enumerate}
Sometimes by abuse of notation, $V$ will be referred to as a virtual
crystal.

Let $b\in \Vh$ and $i\in \I^X$. We say that $b$ is $i$-aligned if
\begin{enumerate}
\item $\vp{j}^Y(b)=\vp{j'}^Y(b)$ for all $j,j'\in\bij(i)$, and similarly
for $\ve{}$.
\item $\vp{j}^Y(b)\in \mult_i \Z$ for all $j\in\bij(i)$ and similarly for $\ve{}$.
\end{enumerate}
In this case
\begin{equation}
  \vp{i}^X(b)=\frac{1}{\mult_i} \vp{j}^Y(\viso(b)) \qquad\text{for $j\in\bij(i), b\in B$}
\end{equation}
and similarly for $\ve{}$. Say that $b\in \Vh$ is aligned if it is $i$-aligned for all
$i\in \I^X$ and a subset $V\subset \Vh$ is aligned if all its elements are.

\begin{prop} \cite{OSS:2001} Aligned virtual crystals form a tensor category.
\end{prop}

Say that $(\V,\Vh)$ is simple if $\V$ and $\Vh$ are simple crystals.
For the rest of the definitions we assume that
the virtual crystals are simple and aligned.

Let $(\V,\Vh)$ and $(\V',\Vh')$ be virtual crystals.
\begin{defconj}\label{defc:Rv}
Define the virtual $R$-matrix $\Rv:\V\otimes \V' \to \V'\otimes V$ as the restriction
of the type $Y$ $R$-matrix $\Rh:\Vh\otimes \Vh' \to \Vh'\otimes \Vh$.
\end{defconj}
For this definition to make sense it needs to be shown that
$\Rh(\V\otimes \V') \subset \V'\otimes \V$. In this case,
let $\viso:B\cong V$ and $\viso':B'\cong V'$
be $X$-crystal isomorphisms. By the uniqueness of the
$R$-matrix it follows that the diagram
\begin{equation} \label{eq:vR}
\begin{CD}
  B \otimes B' @>R>> B' \otimes B \\
  @V{\viso\otimes\viso'}VV @VV{\viso'\otimes\viso}V \\
  V \otimes V' @>{\Rv}>> V' \otimes V
\end{CD}
\end{equation}
commutes.

\begin{definition}
Define the virtual energy function $\Hv:\V\otimes \V'\to \Z$ by
\begin{equation*}
  \Hv(b \otimes b') = \frac{1}{\mult_0} H_Y(b\otimes b')
\end{equation*}
where $H_Y:\Vh\otimes \Vh'\to \Z$.
\end{definition}
If Definition-Conjecture \ref{defc:Rv} holds, it follows that
\begin{equation} \label{eq:vH}
  H_X(b\otimes b') = \Hv(\viso(b)\otimes \viso'(b'))
\end{equation}
where $H_X:B\otimes B'\rightarrow\Z$ is the energy function.

Similarly, define $\Dv:\V\to\Z$ as
\begin{equation*}
\Dv(b)=\frac{1}{\mult_0} D_{\Vh}(b).
\end{equation*}
If \eqref{eq:b nat} and Definition-Conjecture \ref{defc:Rv} hold then
\begin{equation}
  D_X(b) = \Dv(\viso(b)) \qquad\text{for $b\in B$.}
\end{equation}
where $D_X:B\rightarrow\Z$ is the intrinsic energy of $B$.

Finally, let $\la\in\Pfin^+$ for the algebra $X$ and
\begin{equation*}
\Path(\V,\la)=\{b\in \V \mid \text{$\wt(b)=\emb(\la)$ and
$\eh_i b=0$ undefined for $i\in \J^X$}\}.
\end{equation*}
Then let
\begin{equation*}
\Xv(\V,\la)=\sum_{b\in \Path(\V,\la)} q^{\Dv(b)}.
\end{equation*}

Let us define the $Y$-crystal
\begin{equation*}
\Vh^{r,s} = \bigotimes_{j\in \bij(r)} B_Y^{j,\mult_r s}
\end{equation*}
except for $A_{2n}^{(2)}$ and $r=n$ in which case
$\Vh^{n,s} = B_Y^{n,s}\otimes B_Y^{n,s}$.

\begin{definition} Let $\V^{r,s}$ be the subset of
$\Vh^{r,s}$ generated from $\uu(\Vh^{r,s})$ using the virtual
crystal operators $\eh_i$ and $\fh_i$ for $i \in \I^X$.
\end{definition}

\begin{conjecture} \label{conj:crystal}
\mbox{}
\begin{enumerate}
\item[(V1)] The pair $(\V^{r,s},\Vh^{r,s})$ is a simple aligned virtual crystal.
\item[(V2)] There is an isomorphism of $X$-crystals
\begin{equation*}
\viso:B^{r,s}_X \cong \V^{r,s}
\end{equation*}
such that $e_i$ and $f_i$ correspond to $\eh_i$ and $\fh_i$ respectively,
for all $i\in I^X$.
\item[(V3)] Let $\la$ be a classical dominant weight for $X$, $B$
a tensor product of $X$-crystals of the form $B^{r,s}$, and
$(\V,\Vh)$ the corresponding tensor product of virtual crystals $(\V^{r,s},\Vh^{r,s})$.
Then
\begin{equation}
  \Xv(\V,\la)=X(B,\la).
\end{equation}
\end{enumerate}
\end{conjecture}

In \cite{OSS:2001} Conjecture \ref{conj:crystal} is proved for embeddings
$C_n^{(1)}, A_{2n}^{(2)}, A_{2n}^{(2)\dagger}, D_{n+1}^{(2)}
\hookrightarrow A_{2n-1}^{(1)}$ and
tensor factors of the form $B^{r,1}$.

\begin{theorem} \label{thm:row crystal} Conjecture \ref{conj:crystal}
holds when $X$ is of nonexceptional affine type and $B$ is a tensor
product of crystals of the form $B^{1,s}$.
\end{theorem}
This theorem is proven in subsections \ref{subsec:row crystal D}
and \ref{subsec:row crystal A}.

\subsection{Virtual crystals $V^{1,s}$ for
$A_{2n-1}^{(2)},B_n^{(1)}\hookrightarrow D_{n+1}^{(1)}$}
\label{subsec:row crystal D}

\begin{prop} \label{pp:row crystal} For $X=B_n^{(1)}$ and $Y=D_{n+1}^{(1)}$,
\begin{equation*}
V^{1,s} = \{b\in \Vh^{1,s} \mid \text{$x_i,\xb_i\in 2\Z$ for $i<n$,
 $x_n+\xb_n\in 2\Z$, $x_{n+1}=\xb_{n+1}=0$} \}.
\end{equation*}
Moreover Theorem \ref{thm:row crystal} holds.
\end{prop}
\begin{proof} The explicit form of $V^{1,s}$ follows from
$\uu(\Vh^{1,s})=1^{2s}$ and the definitions of the virtual
crystal operators. It is easy to show that for $s=1$ the map
$B^{1,1}\rightarrow \V^{1,1}$ defined by
$i\mapsto ii$ and $\bar{i}\mapsto \bar{i}\bar{i}$ for $1\le i\le n$
and $\0\mapsto n \bar{n}$, is the desired isomorphism for $s=1$.
Similarly, it is straightforward to show that for $s$ arbitrary, the
desired isomorphism $\viso:B^{1,s}\rightarrow \V^{1,s}$ is given by
replacing each letter (which is an element of $B^{1,1}$) of a word in $B^{1,s}$
by the corresponding pair of letters as in the case $s=1$. This proves (V1) and
(V2).

For (V3) we need to check that
\begin{equation}
  \Dv(\viso(b)) = D_B(b)\qquad\text{for $b\in B$.}
\end{equation}
Since $D_B$ is defined in terms of $R$, $H$ and functions
$D_{B^{1,s}}$, it suffices to verify \eqref{eq:b nat} and
Definition-Conjecture \ref{defc:Rv}.

The element $\nat(B^{1,s})$ is given explicitly by $\bar{1}^{s}$.
By the explicit computation of $H:B^{1,s}\otimes B^{1,s}\rightarrow\Z$
given in \cite{HKOT:2002} it follows that \eqref{eq:b nat} holds.

To check Definition-Conjecture \ref{defc:Rv} we consider the explicit
expressions for the $R$-matrices of types $B_n^{(1)}$ and $D_{n+1}^{(1)}$
given in \cite{HKOT:2002}. From this it suffices to show that the images of
relations in the plactic monoid of type $B_n$ are relations in the
plactic monoid of type $D_{n+1}$ \cite{L:2001}. This is straightforward.
\end{proof}

\begin{prop} For $X=A_{2n-1}^{(2)}$ and $Y=D_{n+1}^{(1)}$,
\begin{equation*}
V^{1,s} = \{b\in \Vh^{1,s} \mid \text{$x_{n+1}=\xb_{n+1}=0$} \}.
\end{equation*}
Moreover Theorem \ref{thm:row crystal} holds.
\end{prop}
\begin{proof} The proof is similar to that of Proposition \ref{pp:row crystal}.
In particular the bijection $B^{1,s}\rightarrow V^{1,s}$ is given by leaving a word
unchanged.
\end{proof}

\subsection{Virtual crystals $V^{1,s}$ for
$C_n^{(1)}, A_{2n}^{(2)}, A_{2n}^{(2)\dagger}, D_{n+1}^{(2)}\hookrightarrow
A_{2n-1}^{(1)}$}
\label{subsec:row crystal A}
We require some preliminaries on crystals of type $A_{2n-1}^{(1)}$.

Consider $Y=A_{2n-1}^{(1)}$ and $X$ one of $C_n^{(1)},
A_{2n}^{(2)}, A_{2n}^{(2)\dagger},D_{n+1}^{(2)}$. In all these cases $\Vh^{1,s} =
B_Y^{2n-1,s} \otimes B_Y^{1,s}$. We introduce the alphabets
\begin{equation}
  \Y = \{ 1<2<\dotsm<2n \} \qquad
 \Yv = \{ 2n^\vee < (2n-1)^\vee <\dotsm < 2^\vee < 1^\vee \}.
\end{equation}
$\Y$ and $\Yv$ are the sets of elements of $B_Y^{1,1}$ and
$(B_Y^{1,1})^\vee \cong B_Y^{2n-1,1}$ respectively. The element
$i^\vee\in B_Y^{2n-1,1}$ is the column of height $2n-1$ in the alphabet
$\Y$ with the letter $i$ missing. For $1\le i\le 2n-1$,
$f_i((2n+1-i)^\vee)=(2n-i)^\vee$ and $f_i(b)$ is undefined otherwise.
$f_0(1^\vee)=(2n)^\vee$ and $f_0(b)$ is undefined otherwise.
In this notation, $B_Y^{2n-1,s}$ consists of the weakly increasing words
of length $s$ in the alphabet $\Yv$. For $b=b_1\otimes b_2 \in \Vh^{1,s}$,
let $y_i$ be the number of letters $i$ in $b_2$ and $\yv_i$ the number of
letters $i^\vee$ in $b_1$, for $1\le i\le 2n$.

The $R$-matrix $R:B_Y^{1,1} \otimes B_Y^{2n-1,1} \rightarrow
B_Y^{2n-1,1}\otimes B_Y^{1,1}$ is given by
\begin{equation} \label{eq:comm}
  i \otimes j^\vee \rightarrow
  \begin{cases}
    j^\vee \otimes i & \text{if $i\not=j$} \\
    (i+1)^\vee \otimes (i+1) & \text{if $i=j<2n$} \\
    1^\vee \otimes 1 & \text{if $i=j=2n$.}
  \end{cases}
\end{equation}
The $R$-matrix $R:B_Y^{1,s} \otimes B_Y^{2n-1,s}\rightarrow B_Y^{2n-1,s}
\otimes B_Y^{1,s}$ is given by iterating the above $R$-matrix so that all
of the elements of $\Yv$ are commuted to the left. The element
$1^\vee \otimes 1$ commutes with all elements of $B_Y^{1,1}$
and $B_Y^{2n-1,1}$.

To formulate the next propositions we also need an involution
$*:B\to B$ on crystals of type $A_{2n-1}^{(1)}$ \cite[Section 3.8]{OSS:2001}.
Given a word $u$, let $u^*$ be the word obtained by replacing
each letter $i$ by $2n+1-i$, and reversing the resulting word.
Clearly if $u$ is a column word then so is $u^*$. If $b=c_1
c_2\dots c_s\in B^{r,s}$ where $c_j$ is a column word for all $j$, then by
definition $b^*=c_s^* \dots c_1^*\in B^{r,s}$, which is a sequence of
column words. Under this map the crystal operators transform as follows:
\begin{equation*}
\begin{split}
  f_i(b^*) &= e_{n-i}(b)^* \\
  e_i(b^*) &= f_{n-i}(b)^* \\
  \wt(b^*) &= w_0 \wt(b).
\end{split}
\end{equation*}

\begin{prop} For $X=C_{n}^{(1)}$ and $Y=A_{2n-1}^{(1)}$,
\begin{equation} \label{eq:v C row}
V^{1,s} = \{b\in \Vh^{1,s} \mid \text{$b^{\vee*}=R(b)$, $\min(y_1,\yv_1),\min(y_{n+1},\yv_{n+1})\in 2\Z$}\}
\end{equation}
Moreover Theorem \ref{thm:row crystal} holds.
\end{prop}
\begin{proof} We first prove \eqref{eq:v C row}. By the definition of $V^{1,s}$,
it suffices to show that the right hand side $V'$ of \eqref{eq:v C row}
contains $\uu(\Vh^{1,s})$, and every element of $V'$ is reachable from $\uu(\Vh^{1,s})$
using the virtual crystal operators $\eh_i,\fh_i$ for $i\in \I^X$.

We first digress on the self-duality condition
\begin{equation} \label{eq:selfdual}
  b^{\vee*}=R(b).
\end{equation}
By the proof of \cite[Prop 6.8]{OSS:2001}, in the set $\Vh^{r,s}$, the condition
\eqref{eq:selfdual}
is preserved under $e_0$, $e_n$, $\eh_i$ for $1\le i\le n-1$, and similarly for $f$.
For $b\in \Vh^{1,s}$, using \eqref{eq:comm} equation \eqref{eq:selfdual} is equivalent to
\begin{align}
\label{eq:dual1}  y_{2n+1-i} &= \yv_i - \min(y_i,\yv_i) + \min(y_{i+1},\yv_{i+1}) \\
\label{eq:dual2} \yv_{2n+1-i} &= y_i - \min(y_i,\yv_i) + \min(y_{i+1},\yv_{i+1})
\end{align}
for $1\le i\le 2n$, where $y_{2n+1}=y_1$ and $\yv_{2n+1}=\yv_1$.

We deduce two consequences of \eqref{eq:selfdual}.
Subtracting \eqref{eq:dual1} and \eqref{eq:dual2} we obtain
\begin{equation} \label{eq:sum comp}
  y_i + y_{2n+1-i} = \yv_i + \yv_{2n+1-i}
\end{equation}
for $1\le i\le 2n$. We also have $\ve{0}(b)=y_1+\yv_{2n}-\min(y_{2n},\yv_{2n})$.
By \eqref{eq:dual2} with $i=2n$ and \eqref{eq:sum comp} with $i=1$, we have
\begin{equation} \label{eq:even0}
  \ve{0}(b)=2 y_1 - \min(y_1,\yv_1).
\end{equation}

Now we show that $u=\uu(\Vh^{1,s})\in V'$. This element satisfies
$y_i = s \delta_{i,1}$ and $\yv_i = s \delta_{i,2n}$ for $1\le i\le 2n$.
Comparing this with \eqref{eq:dual1} and \eqref{eq:dual2} it follows that
$u$ satisfies \eqref{eq:selfdual}. It follows that $u\in V'$.

We next check that $V'$ aligned. Let $b\in V'$ and $i\in \I^X$.
Since $b$ satisfies \eqref{eq:selfdual} it is $i$-aligned if $1\le
i\le n-1$ by \cite[Prop. 6.9]{OSS:2001}. For $0$-alignedness, by
\eqref{eq:even0} we see that $\ve{0}(b)$ is even since
$\min(y_1,\yv_1)$ is. The proof that $\vp{0}(b)$ is even is
similar. So $b$ is $0$-aligned. The proof that $b$ is $n$-aligned,
is similar as well. So $V'$ is aligned.

Next it is shown that the set $V'$ is closed under
$\eh_i$ and $\fh_i$ for $i\in \I^X$. Let $b\in V'$.
$\eh_i b$ is self-dual since $b$ is. Note that the quantity $\min(y_1,\yv_1)$
is unchanged for $i\not\in\{0,1\}$.
We have $\ve{1}(b_1)=\yv_1$ and $\vp{1}(b_2)=y_1$.
Hence by the tensor product rule,
$\min(y_1,\yv_1)$ remains the same upon applying $\eh_1$.
Let $i=0$. Since $b\in V'$, $b$ is $0$-aligned, so that
$\ve{0}(b)\in 2\Z$. Since $\ve{0}(\eh_0 b)=\ve{0}(b)-2$ is even,
by \eqref{eq:even0}, the self-dual element $\eh_0 b$ has the
property that $\min(y_1,\yv_1)\in 2\Z$. Thus $\eh_i b$
satisfies that property for all $i$. The property that
$\min(y_{n+1},\yv_{n+1})\in2\Z$ is satisfied for $\eh_i b$
is similar. Thus $\eh_i b\in V'$ for all $i\in \I^X$.
The proof that $\fh_i b\in V'$ for all $i\in \I^X$ is again similar.

Let $b\in V'$. It suffices to find a sequence of operators $\eh_i$
and $\fh_i$ leading from $b$ to $u$. We shall induct on the
quantity $\min(y_1,\yv_1)$, which is invariant under $\eh_i$ and
$\fh_i$ for $i\in \I^X\backslash\{0\}$ by previous arguments.
Suppose first that $\ve{j'}(b)>0$ for some $j'\not=0$. By
alignedness it follows that we may apply a sequence of operators
$\eh_i$ for $i\in \I^X\backslash\{0\}$ to $b$, thereby passing to a
classical highest weight vector of $\Vh^{1,s}$. The classical
highest weight vectors of $\Vh^{1,s}$ are given explicitly by $u_k
= (2n^\vee)^{s-k} 1^{\vee k} \otimes 1^s$, for $0\le k\le s$.
$u_k$ satisfies $\min(y_1,\yv_1)=k$. By assumption $b=u_k$ for $k$
even. If $k=0$ then $b=u_0=u$ and we are done. If $k>0$ then
$\fh_0 b$ satisfies $\min(y_1,\yv_1)=k-2$, which is even. We are
done by induction.

We have shown that \eqref{eq:v C row} holds and that $V^{1,s}$ is aligned.

The bijection $\viso:B^{1,s}\rightarrow V^{1,s}$ is given as
follows. Let $b\in B^{1,s}$.
In the case $s=1$, the map $B^{1,1}\rightarrow V^{1,1}$ is given
by $i \mapsto (2n+1-i)^\vee \otimes i$ and
$\bar{i}\mapsto i^\vee \otimes (2n+1-i)$. The map $\viso:B^{1,s}\rightarrow V^{1,s}$
is given by the composite map
\begin{equation} \label{eq:crystal map}
B^{1,s} \hookrightarrow
  (B^{1,1})^{\otimes s} \rightarrow (B_Y^{2n-1,1} \otimes B_Y^{1,1})^{\otimes s}
  \rightarrow (B_Y^{2n-1,1})^{\otimes s} \otimes (B_Y^{1,1})^{\otimes s}.
\end{equation}
It follows from \eqref{eq:comm} that the image of this map is contained in
$B_Y^{2n-1,s} \otimes B_Y^{1,s}$. Computing this commutation
explicitly and using the notation $x_i,\xb_i$ to
describe $b$ for $1\le i\le n$, and $y_i,\yv_i$ for
$\viso(b)$, we have
\begin{equation} \label{eq:xy embed}
\begin{split}
  y_1 &= x_1 - \min(x_1,\xb_1) + s-\sum_{i=1}^n
  (x_i+\xb_i) \\
  \yv_1 &= \xb_1 - \min(x_1,\xb_1) + s-\sum_{i=1}^n
  (x_i+\xb_i) \\
  y_i &= x_i - \min(x_i,\xb_i) +
  \min(x_{i-1},\xb_{i-1})\qquad\text{for $i>1$} \\
  \yv_i &= \xb_i - \min(x_i,\xb_i) +
  \min(x_{i-1},\xb_{i-1})\qquad\text{for $i>1$}
\end{split}
\end{equation}
To recover $y_i$ and $\yv_i$ for $n+1\le i\le 2n$ one may use
\eqref{eq:dual1} and \eqref{eq:dual2}, plus the fact that the
total number of letters in either $b_1$ or $b_2$, is $s$.

The composite map given in \eqref{eq:crystal map} sends $e_i$ to $\eh_i$
for $1\le i\le n$ \cite{B:2000}. It is straightforward to check that
$e_0$ goes to $\eh_0$ using
\eqref{eq:xy embed}, \eqref{eq:dual1}, and \eqref{eq:dual2}.
Therefore $\viso$ is a morphism of $X$-crystals. It is clearly injective.
The image is $V^{1,s}$ since $\viso(\uu(B^{1,s}))=\uu(\Vh^{1,s})$ and both $B^{1,s}$ and
$V^{1,s}$ are connected. Therefore $\viso:B^{1,s}\rightarrow V^{1,s}$ is an isomorphism
of $X$-crystals. This completes the proof of (V1) and (V2).
(V3) follows by \cite[Section 6.6]{OSS:2001}.
\end{proof}

\begin{prop} For $X=A_{2n}^{(2)}$ and $Y=A_{2n-1}^{(1)}$,
\begin{equation} \label{eq:v A2 row}
V^{1,s} = \{b\in \Vh^{1,s} \mid
 \text{$b^{\vee*}=R(b)$, $\min(y_{n+1},\yv_{n+1})\in 2\Z$}\}.
\end{equation}
Moreover Theorem \ref{thm:row crystal} holds.
\end{prop}

The proof is entirely similar to that of $C_n^{(1)}$.

\begin{prop} For $X=D_{n+1}^{(2)}$ and $Y=A_{2n-1}^{(1)}$,
\begin{equation} \label{eq:v D2 row}
V^{1,s} = \{b\in \Vh^{1,s} \mid \text{$b^{\vee*}=R(b)$}\}.
\end{equation}
Moreover Theorem \ref{thm:row crystal} holds.
\end{prop}
\begin{proof} For $X=D_{n+1}^{(2)}$ most of the proof is similar to that of type $C_n^{(1)}$.
Here the classical subalgebra of $X$ is of type $B_n$, so the isomorphism
$\viso:B^{1,s}\rightarrow V^{1,s}$ is a bit different. It is given by
$\0\mapsto (n+1)^\vee \otimes (n+1)$, with the other letters mapped as in the $C_n^{(1)}$ case.
The explicit map is given as in \eqref{eq:xy embed} except that
\begin{equation} \label{eq:xy embed D}
\begin{split}
  y_1 &= x_1 - \min(x_1,\xb_1) + s-x_\0-\sum_{i=1}^n
  (x_i+\xb_i) \\
  \yv_1 &= \xb_1 - \min(x_1,\xb_1) + s-x_\0-\sum_{i=1}^n
  (x_i+\xb_i).
\end{split}
\end{equation}
\end{proof}

\begin{prop} For $X=A_{2n}^{(2)\dagger}$ and $Y=A_{2n-1}^{(1)}$,
\begin{equation} \label{eq:v A2d row}
V^{1,s} = \{b\in \Vh^{1,s} \mid \text{$b^{\vee*}=R(b),\min(y_1,\yv_1)\in 2\Z$}\}.
\end{equation}
Moreover Theorem \ref{thm:row crystal} holds.
\end{prop}
The proof is similar.

\section{Fermionic formula}
\label{sec:fermionic}

\subsection{Review}
\label{subsec:Mreview} This subsection reviews definitions of
\cite{HKOTT:2002,HKOTY:1999}. For this section we assume that
$\geh\not=A_{2n}^{(2)\dagger}$; for that type we refer the reader
to~\cite[Section 7.6]{OSS:2001}.
Fix $\la\in\Pfin^+$ and $B$ a tensor product of crystals of the form $B^{r,s}$.
Let $L_i^{(a)}$ be the number of tensor factors in $B$ that are equal to
$B^{a,i}$. Set $\alt_a=\alpha_a$ for all $a\in \J$ except for type
$A_{2n}^{(2)}$ in which case $\alt_a$ are the simple roots of type
$B_n$.

Let $\nu=(m_i^{(a)})$ be a matrix of nonnegative integers for
$i\in\Z_{>0}$ and $a\in \J$.
Say that $\nu$ is a $(B,\la)$-configuration if
\begin{equation}
\label{eq:config}
\sum_{\substack{a\in \J \\ i\in\Z_{>0}}} i\, m_i^{(a)} \alt_a
= \sum_{\substack{a\in\J \\ i\in\Z_{>0}}} i\, L_i^{(a)} \Lab_a - \la
\end{equation}
except for type $A_{2n}^{(2)}$. In this case the right hand side should
be replaced by $\iota(\mathrm{r.h.s})$ where $\iota$ is a $\Z$-linear
map from the weight lattice of type $C_n$ to the weight lattice of
type $B_n$ such that
\begin{equation*}
\iota(\Lab_a^C)=\begin{cases} \Lab_a^B & \text{for $1\le a<n$}\\
2\Lab_a^B & \text{for $a=n$.}\end{cases}
\end{equation*}

Say that a configuration $\nu$ is admissible if
\begin{equation} \label{eq:ppos}
  p_i^{(a)} \ge 0\qquad\text{for all $a\in\J$ and
  $i\in\Z_{>0}$,}
\end{equation}
where
\begin{equation} \label{eq:p}
p_i^{(a)} = \sum_{k\in\Z_{>0}} \left( L_k^{(a)} \min(i,k) -
\dfrac{1}{t_a^\vee} \sum_{b\in\J} (\alt_a|\alt_b)\min(t_b i,t_a
k)\, m_k^{(b)}\right).
\end{equation}
Write $C(B,\la)$ for the set of admissible
$(B,\la)$-configurations. Define
\begin{equation} \label{eq:cc}
cc(\nu) = \dfrac{1}{2} \sum_{a,b\in\J} \sum_{j,k\in\Z_{>0}} (\alt_a|\alt_b)
\min(t_b j, t_a k) m_j^{(a)} m_k^{(b)}.
\end{equation}
The fermionic formula is defined by
\begin{equation}\label{eq:fermi}
M(B,\la;q) = \sum_{\nu\in C(B,\la)} q^{cc(\nu)}
\prod_{a\in\J} \prod_{i\in\Z_{>0}}
\qbin{p_i^{(a)}+m_i^{(a)}}{m_i^{(a)}}_{q^{t^\vee_a}}.
\end{equation}

The $X=M$ conjecture of \cite{HKOTT:2002,HKOTY:1999} states that
\begin{equation}\label{eq:X=M}
  X(B,\la;q^{-1})=M(B,\la;q).
\end{equation}

The fermionic formula $M(B,\la)$ can be interpreted using
combinatorial objects called rigged configurations.
Denote by $(\nu,J)$ a pair where $\nu=(m_i^{(a)})$ is a matrix and
$J=(J^{(a,i)})$ is a matrix of partitions with $a\in\J$ and $i\in\Z_{>0}$.
Then a rigged configuration is a pair $(\nu,J)$ such that $\nu\in C(B,\la)$
and the partition $J^{(a,i)}$ is contained in a $m_i^{(a)}(\nu) \times
p_i^{(a)}(\nu)$ rectangle for all $a,i$. The set of rigged
$(B,\la)$-configurations for
fixed $\la$ and $B$ is denoted by $\RC(B,\la)$. Then
\eqref{eq:fermi} is equivalent to
\begin{equation*}
M(B,\la)=\sum_{(\nu,J)\in\RC(B,\la)} q^{\cc(\nu,J)}
\end{equation*}
where $\cc(\nu,J)=\cc(\nu)+|J|$
and $|J|=\sum_{(a,i)} t_a^\vee |J^{(a,i)}|$.

\subsection{Virtual fermionic formula}

We define virtual rigged configurations in analogy to virtual
crystals.

\begin{definition} \label{def:VRC}
Let $X$ and $Y$ be as in \eqref{eq:embed}, and $\la$ and $B$
as in subsection \ref{subsec:Mreview} for type $X$. Let $(V,\Vh)$ be the virtual
$Y$-crystal corresponding to $B$.
Then $\RCv(B,\la)$ is the set of elements $(\nh,\Jh)\in \RC(\Vh,\emb(\la))$
such that:
\begin{enumerate}
\item For all $i\in \Z_{>0}$,
$\mh_i^{(a)}=\mh_i^{(b)}$ and $\Jh_i^{(a)}=\Jh_i^{(b)}$ if $a$ and $b$ are in the same $\aut$-orbit in $I^Y$.
\item For all $i\in \Z_{>0}$, $a\in \J^X$, and $b\in \bij(a)\subset I^Y$, we have
$\mh_j^{(b)}=0$ if $j \not\in \mult_a \Z$ and the parts of $\Jh_i^{(b)}$ are multiples
of $\mult_a$.
\end{enumerate}
\end{definition}

\begin{theorem} \label{thm:M=VM} There is a bijection $\RC(B,\la)\rightarrow \RCv(B,\la)$
sending $(\nu,J)\mapsto(\nh,\Jh)$ given
as follows. For all $a\in \J^X$, $b\in\bij(a)\subset \J^Y$, and
$i\in\Z_{>0}$,
\begin{align}
\mh_{\mult_a i}^{(b)} &= m_i^{(a)} \\
\Jh_{\mult_a i}^{(b)}&=\mult_a J_i^{(a)},
\end{align}
except when $X=A_{2n}^{(2)}$ and $a=n$, in which case
\begin{equation*}
\begin{split}
  \mh_i^{(n)} &= m_i^{(n)} \\
  \Jh_i^{(n)} &= 2 J_i^{(n)}.
\end{split}
\end{equation*}
The cocharge changes by
\begin{equation} \label{eq:ccembed}
\cc(\nh,\Jh) = \mult_0 \,\cc(\nu,J).
\end{equation}
\end{theorem}
\begin{proof} Let $\Lh$ be to $\Vh$ as $L$ is to $B$ as in
subsection \ref{subsec:Mreview}.
For $a\in \J^X$, $b\in\bij(a)$, and $i\in\Z_{>0}$,
\begin{equation*}
\begin{split}
  \Lh_{\mult_a i}^{(b)} &= L_i^{(a)} \\
  \Lh_j^{(b)} &= 0 \qquad\text{for $j\not\in \mult_a \Z$},
\end{split}
\end{equation*}
except when $X=A_{2n}^{(2)}$ and $a=n$, in which case
$\Lh_i^{(n)}=2 L_i^{(n)}$ for all $i$.
Using \eqref{eq:p} we have, for all $b\in \bij(a)$ and $i\in\Z_{>0}$,
\begin{equation*}
  \ph_{\mult_a i}^{(b)} = \mult_a p_i^{(a)},
\end{equation*}
except when $X=A_{2n}^{(2)}$ and $i=n$, in which case
$\ph_i^{(n)}=2 p_i^{(n)}$.
Therefore $(\nu,J)\mapsto(\nh,\Jh)$ defines a bijection. Using
\eqref{eq:cc} we see that \eqref{eq:ccembed} holds.
\end{proof}

\section{Algorithms for computing the fermionic formula}
\label{sec:kleber}

To compute the fermionic formula $M(B,\la)$, one must find the set of
admissible $(B,\la)$-configurations $C(B,\la)$. One direct approach would be
to test the admissibility conditions \eqref{eq:ppos}
on the set of $(B,\la)$-configurations \eqref{eq:config}
which consist of all possible $n$-tuples of partitions of sizes that
depend on $\la$ and $B$. This quickly becomes infeasible as $B$ and $\la$
grow.

In \cite{Kl:1997,Kl:1998} Kleber gives an efficient algorithm to compute
the set of admissible configurations in the simply-laced types
$A_n^{(1)}$, $D_n^{(1)}$, $E_6^{(1)}$, $E_7^{(1)}$, and $E_8^{(1)}$.
It generates a rooted tree $T(B)$ whose nodes are labelled by elements
of $\Pfin^+$. The tree $T(B)$ is constructed to have the property that
the elements of $C(B,\la)$ are in bijection with the nodes of $T(B)$ labelled $\la$.
If a node $x$ labelled $\la$ corresponds to a configuration $\nu$,
then $\nu$ can be recovered from the unique path in $T(B)$ from $x$ to the root.

\subsection{Kleber's algorithm}
We review Kleber's algorithm \cite{Kl:1997,Kl:1998}.
Let $X$ be the Dynkin diagram of an untwisted affine Lie algebra
whose canonical simple subalgebra is of simply-laced type.
Let $B$ and $L$ be as in subsection \ref{subsec:Mreview}.

We define a tree $T(B)$ by the following algorithm. Each node $x$ is labelled by
an element $\wt(x)\in \Pfin^+$ called its weight. It has the property that if $x$ is a node
and $y$ is its child, then $\wt(x)\not=\wt(y)$ and $\wt(x)\dom\wt(y)$. A tree edge $(x,y)$
is labelled by the element $\dd_{xy}=\wt(x)-\wt(y)\in \Qfin^+\backslash \{0\}$.
\begin{enumerate}
\item Let $T_0$ be the tree consisting of a single node of weight $0$
and set $\ell=0$.
\item \label{it:inc} Add 1 to $\ell$.
\item \label{it:addcols} Let $T_\ell'$ be obtained from $T_{\ell-1}$ by adding
$\sum_{a=1}^n \Lab_a \sum_{i\ge \ell} L_i^{(a)}$ to the weight of each node.
\item \label{it:movedown} Let $T_\ell$ be obtained from $T_\ell'$ as follows.
Let $x$ be a node at depth $\ell-1$ of weight $\mu$.
Suppose there is a weight $\tau\in\Pfin^+$ such that $\mu\not=\tau$,
$\mu\dom\tau$, and if $x$ is not the root,
$\nu-2\mu+\tau\in \Qfin^+$ where $\nu$ is the weight
of the parent $w$ of $x$.
In every such case we attach to $x$ a child $y$ of weight $\tau$. Note that
if $x$ is not the root, the condition $\nu-2\mu+\tau\in \Qfin^+$
is equivalent to $d_{wx}\dom d_{xy}$.
\item If $T_\ell\not=T_{\ell-1}$ then go to step \ref{it:inc}.
\item Otherwise set $T(B)=T_\ell$ and stop.
\end{enumerate}
For large $\ell$ step \ref{it:addcols} does not change the tree. For such
$\ell$, step \ref{it:movedown} can only be applied finitely many times
since there are finitely many elements of $\Pfin^+$ dominated by a given
element of $\Pfin^+$. Hence the algorithm terminates.

There is a bijection from the nodes of $T(B)$ and the
configurations $C(B)=\bigcup_{\la\in\Pfin^+} C(B,\la)$ given
as follows. Let $x$ be a node at depth $p$ in $T(B)$ of weight $\la$.
Let $\la^{(0)},\la^{(1)},\dotsc,\la^{(p)}=\la$ be the weights of the
nodes on the path from the root of $T(B)$ to $x$.
Then the configuration $\nu\in C(B,\la)$ corresponding to $x$ is defined by
\begin{equation}
  m_i^{(a)} = (\la^{(i-1)}-2\la^{(i)}+\la^{(i+1)}\mid \Lab_a)
\end{equation}
where we make the convention that $\la=\la^{(p+1)}=\la^{(p+2)}=\dotsm$.
The vacancy numbers are given by
\begin{equation}\label{eq:vac}
  p_i^{(a)} = -\sum_{j>i} (j-i) L_j^{(a)} + (\la^{(i)}\mid\alpha_a).
\end{equation}

Suppose we are only interested in finding $C(B,\la)$ for a particular
$\la\in\Pfin^+$. It is wasteful to generate the entire tree $T(B)$
and then select the nodes of weight $\la$. Because the weight of a node
dominates that of any of its children,
we can prune the tree as follows. In step \ref{it:movedown}, we only add a node
of weight $\tau$ at depth $\ell$ if
\begin{equation} \label{eq:trimdom}
  \tau':=\tau + \sum_{j>\ell} \Lab_a L_j^{(a)} \dom \la.
\end{equation}
There is another condition under which we can prune.
Suppose that in the absence of pruning, we would have added
a node $y$ of weight $\tau$ at depth $\ell$ in step \ref{it:movedown},
with parent $x$. Then we do not add $y$
if there is an $a$ such that
$(\tau'-\la\mid\Lab_a)>0$ and $(d_{xy} \mid \Lab_a)=0$.
For in this case, the condition in step \ref{it:movedown}
prevents one from reaching the weight $\la$ as a descendant of $\tau$.

\begin{example} Let $B=B^{3,2}\otimes B^{2,1}\otimes B^{1,1}\otimes B^{1,1}$
of type $A_3^{(1)}$. The Kleber algorithm produces the tree $T(B)$
given in Figure \ref{fig:kleber}.
The corresponding configurations are given in the following diagram, where we
represent $\nu$ as a sequence of partitions $\nu^{(a)}$ with $m_i^{(a)}$
rows of length $i$. The vacancy number $p_i^{(a)}$ is placed to the right
of a row of length $i$ in $\nu^{(a)}$.

\Yboxdim{9pt}
\Yinterspace{1pt}
\begin{equation*}
\begin{array}{c|lll}
\la & \nu^{(1)} & \nu^{(2)} & \nu^{(3)}\\ \hline &&&\\
2\Lab_1+\Lab_2+2\Lab_3 & \es & \es & \es\\[3mm]
2\Lab_2+2\Lab_3 & \yngrc(1,0) & \es & \es \\[3mm]
\Lab_1+3\Lab_3 & \yngrc(1,1) & \yngrc(1,0) & \es \\[3mm]
\Lab_1+\Lab_2+\Lab_3 & \yngrc(1,1) & \yngrc(1,1) & \yngrc(1,0) \\[3mm]
2\Lab_3
 & \yngrc(1,0,1,) & \yngrc(1,0,1,) & \yngrc(1,1) \\[5mm]
 & \yngrc(2,0) & \yngrc(2,0) & \yngrc(1,0) \\[3mm]
3\Lab_1+\Lab_3 & \es & \yngrc(1,0) & \yngrc(1,0) \\[3mm]
\Lab_2
 & \yngrc(2,0) & \yngrc(2,1) & \yngrc(2,0) \\[3mm]
 & \yngrc(1,0,1,) & \yngrc(1,0,1,) & \yngrc(2,0) \\[5mm]
2\Lab_1 & \yngrc(1,1) & \yngrc(2,0) & \yngrc(2,0)
\end{array}
\end{equation*}
\Yboxdim{5pt}
\Yinterspace{1pt}
\begin{figure}
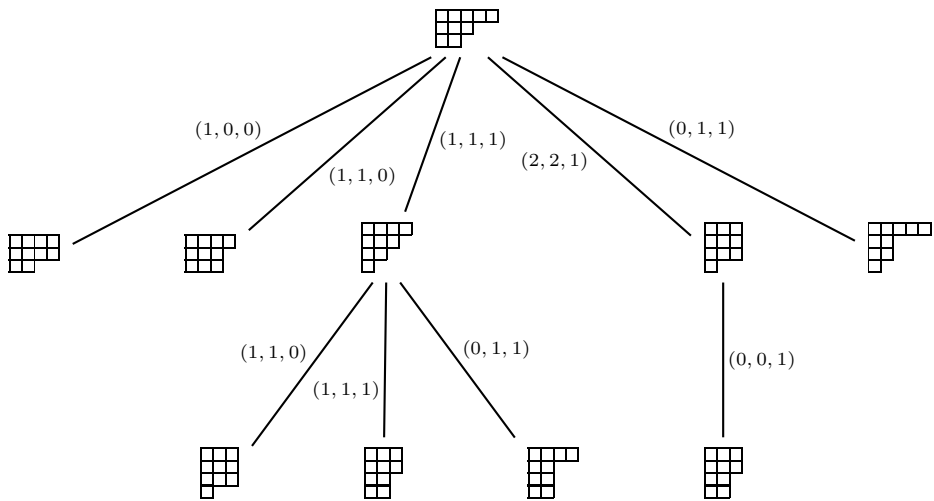

\begin{center}
\psset{levelsep=3cm,nodesep=4pt,treesep=1.6cm,labelsep=1pt,npos=0.5}
\pstree{\TR{\yng(5,3,2)}}
{
\TR{\yng(4,4,2)} \nbput{\scriptsize{$(1,0,0)$}}
\TR{\yng(4,3,3)} \naput[npos=.6]{\scriptsize{$(1,1,0)$}}
\pstree{\TR{\yng(4,3,2,1)} \naput{\scriptsize{$(1,1,1)$}}}
 {
 \TR{\yng(3,3,3,1)} \nbput{\scriptsize{$(1,1,0)$}}
 \TR{\yng(3,3,2,2)} \nbput[npos=.7]{\scriptsize{$(1,1,1)$}}
 \TR{\yng(4,2,2,2)} \naput{\scriptsize{$(0,1,1)$}}
 }
\pstree{\TR{\yng(3,3,3,1)} \nbput{\scriptsize{$(2,2,1)$}}}
 {
 \TR{\yng(3,3,2,2)} \naput{\scriptsize{$(0,0,1)$}}
 }
\TR{\yng(5,2,2,1)} \naput{\scriptsize{$(0,1,1)$}}
}
\end{center}
\caption{\label{fig:kleber} Tree $T(B)$}
\end{figure}
\end{example}

\subsection{Virtual Kleber algorithm}
Outside of the simply-laced case,
Kleber's algorithm does not directly apply.
However we use the embeddings of affine algebras into those of
simply-laced type, where Kleber's algorithm can be applied.
We call our method the virtual Kleber algorithm.
Let $X$ and $Y$ be as in \eqref{eq:embed}. Theorem \ref{thm:M=VM}
defines a bijection $C(B,\la)\cong \Cv(B,\la)$ where $\Cv(B,\la)$
consists of the configurations $\nh\in C(\Vh,\emb(\la))$ constrained as in
Definition \ref{def:VRC}, or equivalently, the $\nh$ such that
$(\nh,\Jh)\in\RCv(B,\la)$ for some $\Jh$. A naive approach would be to
run Kleber's algorithm to compute the set $C(\Vh,\emb(\la))$
and then to select the desired subset $\Cv(B,\la)$. A more efficient way is to
prune the branches that cannot contain elements of $\Cv(B,\la)$.
This results in a good algorithm
to find $\Cv(B,\la)$ and therefore $M(B,\la)$ for any affine type.

More precisely, one only adds the child $y$ to the node $x$
in step \ref{it:movedown} at depth $\ell$ if:
\begin{enumerate}
\item $(\wt(y)\mid \alpha_a)=(\wt(y) \mid \alpha_b)$ if $a$ and $b$ are
in the same $\aut$-orbit of $I^Y$.
\item If $\ell-1\not\in \mult_a\Z$, then $d_{wx}=d_{xy}$ where $w$ is the
parent of $x$.
\end{enumerate}
These conditions are equivalent to those in Definition \ref{def:VRC}.

Let $\Th(B)$ be the resulting tree. Let $\mult=\max_a \mult_a$.
Then there is a bijection between $\Cv(B,\la)$, and the set of nodes $y$
of weight $\la$ in $\Th(B)$ that satisfy either
of the following conditions:
\begin{enumerate}
\item $y$ is at depth $\ell$ with $\ell\in\mult \Z$, or
\item $(d_{xy} \mid \Lab_a)=0$ for every $a$ such that $1<\mult=\mult_a$,
where $x$ is the parent of $y$.
\end{enumerate}

Observe that for $\ell\not\in \mult\Z$,
there may be nodes at depth $\ell$ in $T_\ell$ whose
weights are not in the image of the embedding $P^X\hookrightarrow P^Y$,
but rather in a superlattice of index $\mult$.
These weights, which cannot appear in the final tree, are necessary
as they allow the virtual Kleber algorithm to reach all of the desired weights.

\begin{example} Let $X=C_2^{(1)}$, $Y=A_3^{(1)}$, $B=B^{1,2}\otimes
B^{1,1}\otimes B^{2,1}$. The virtual Kleber algorithm produces the tree
$\Th(B)$ given in Figure \ref{fig:tree}.
The nodes corresponding to elements of $\Cv(B,\la)$ are circled. We list
the configurations corresponding to the
circled nodes, ordered by increasing depth and then from left to right.
Here we represent $\nu$
as a sequence of partitions $\nu^{(a)}$ with $m_i^{(a)}$ rows of length $i$.
The vacancy number $p_i^{(a)}$
is placed to the right of a row of length $i$ in $\nu^{(a)}$.
\Yboxdim{9pt}
\Yinterspace{1pt}
\begin{equation*}
\begin{array}{c|lll}
  \la & \nu^{(1)} & \nu^{(2)} & \nu^{(3)} \\ \hline&&&\\
  3\Lab_1+\Lab_2 & \es & \es & \es \\[3mm]
  \Lab_1+2\Lab_2 & \yngrc(1,0) & \es & \yngrc(1,0) \\[3mm]
  3\Lab_1  & \yngrc(1,1) & \yngrc(2,0) & \yngrc(1,1) \\[3mm]
  \Lab_1+\Lab_2 &  \yngrc(2,1) & \yngrc(2,2) & \yngrc(2,1) \\[3mm]
  \Lab_1 & \yngrc(2,1,1,1) & \yngrc(2,0,2,) & \yngrc(2,1,1,1) \\[5mm]
  & \yngrc(3,0) & \yngrc(4,0) & \yngrc(3,0)
\end{array}
\end{equation*}
\Yboxdim{5pt}
\Yinterspace{1pt}
\begin{figure}
\begin{center}
\psset{levelsep=2cm,nodesep=4pt,treesep=1cm,npos=0.5,labelsep=1pt}
\pstree{\Toval{\yng(8,5,3)}}
{
\pstree{\TR{\yng(7,5,3,1)}
\nbput{\scriptsize{$(1,1,1)$}}
}
{
\Toval{\yng(7,4,4,1)}
\nbput{\scriptsize{$(0,1,0)$}}
\pstree{\Toval{\yng(6,5,3,2)}
\naput{\scriptsize{$(1,1,1)$}}}
{
\pstree{\TR{\yng(5,5,3,3)}
\naput{\scriptsize{$(1,1,1)$}}
}
{
\Toval{\yng(5,4,4,3)}
\naput{\scriptsize{$(0,1,0)$}}
}
}
}
\pstree{\TR{\yng(6,5,3,2)}
\nbput{\scriptsize{$(2,2,2)$}}
}
{
\Toval{\yng(5,4,4,3)}
\nbput{\scriptsize{$(1,2,1)$}}
}
\Toval{\yng(7,6,2,1)}
\naput{\scriptsize{$(1,0,1)$}}
}
\end{center}
\caption{\label{fig:tree} Tree $\Th(B)$}
\end{figure}
\end{example}

\end{document}